\documentclass{article}
\usepackage{amssymb}
\usepackage[dvips]{graphicx}
\usepackage{latexsym}
\usepackage{geometry}
 \def\qed{\hfill $\Box$}
\newtheorem{lem}{Lemma}
\newtheorem{thm}{Theorem}

\title{Genus one fibered knots in 3-manifolds with reducible genus two Heegaard splittings}
\author{Nozomu Sekino}
\date{}
\begin{document}
\maketitle

\begin{abstract}
We give a necessary and sufficient condition for a simple closed curve on the boundary of a genus two handlebody to decompose the handlebody into $T\times I$ ($T$ is a torus with one boundary component). 
We use this condition to decide whether a simple closed curve on a genus two Heegaard surface is a GOF-knot (genus one fibered knot) which induces the Heegaard splitting. 
By using this, we determine the number and the positions with respect to the Heegaard splittings of GOF-knots in the 3-manifolds with reducible genus two Heegaard splittings. 
This is another proof of results of Morimoto \cite{12} and Baker \cite{2}, \cite{3}.
\end{abstract}

\section{Introduction} \label{sec1}
In \cite{1}, Alexander proved that every closed orientable 3-manifold has a fibered link 
and by Myers \cite{13} and Gonz\'alez-Acu\~na this result was improved so that every closed orientable 3-manifold has a fibered knot. 
Thus focusing on fibered knots is not restricting for a study of 3-manifolds and many works for 3-manifolds have been done 
by using the connection with open-book decompositions, contact structures and so on. 
In this paper, we handle fibered knots whose fiber is a genus one surface. 
Though not every manifold has such knots, they can be relatively easily studied because of their low genus and they might be test cases for higher genus.

A fibered knot in a closed orientable 3-manifold $M$ whose fiber is a torus with one boundary component is called a GOF-knot (genus one fibered knot) in $M$. 
Classically, it is known that all GOF-knot in $S^3$ are the (left and right hand) trefoil and the figure eight knot \cite{5}. 
In \cite{12}, Morimoto investigated the number of GOF-knots in some lens spaces by using their monodromies. 
In \cite{3}, Baker counted GOF-knots in each lens space by using the correspondence between GOF-knots in all 3-manifolds and closed 3-braids with axes in $S^3$. 
He also counted GOF-knots in non-prime 3-manifolds by the same way in his other works \cite{2}. 
By the correspondence, the monodromy of a GOF-knot and its fiber can be calculated, 
however its position in the underlying 3-manifold seems not to be easily found in the context of Heegaard splittings, 
which give a fundamental method for representing a 3-manifold.\\ 
\ The purpose of this paper is to reveal the positions of GOF-knots in some special 3-manifolds by putting them on the (almost unique) standard Heegaard surface.  
We first prepare some terminologies and techniques in Section \ref{sec2} 
and then we prove the relation between GOF-knots and simple closed curves on Heegaard surfaces in Section \ref{sec3}. 
Using the results obtained in Section \ref{sec3}, we investigate GOF-knots in individual cases in Section \ref{sec4}. 
In our method, the position of a GOF-knot is easily found since we put it on a genus two Heegaard surface.
The operation used in Section \ref{sec4} is similar to that in Cho and Koda \cite{7}, \cite{8}, \cite{9}. \\\\

\subsection*{Acknowledgements}
I would like to thank my supervisor, Takuya Sakasai 
for giving me many suggestions for this paper. 
Also, I would like to thank everyone in the laboratory I belong to for providing me with meaningful discussions.\\\\\\

\section{Preliminary} \label{sec2}
\subsection*{Heegaard splittings}
{\it Heegaard splitting } is a method for decomposing a closed orientable 3-manifold into two handlebodies of the same genus. 
The closed orientable surface which is the common boundary of two handlebodies is called a {\it Heegaard surface}. 
The genus of the Heegaard surface of a Heegaard splitting is called the genus of the Heegaard splitting. 
It is known that every closed orientable 3-manifold admits a Heegaard splitting. 
A Heegaard splitting is called {\it reducible} if there is an essential simple closed curve (called a {\it reducing curve}) on the Heegaard surface 
which bounds a disk in each of two handlebodies. 
The sphere in the manifold obtained by pasting two disks bounded by a reducing curve along their boundaries is called a {\it Haken sphere} of the Heegaard splitting. 
If the genus of a reducible Heegaard splitting is greater than one, there must be a separating reducing curve. 
In particular, every reducible genus two Heegaard splitting is decomposed into two genus one Heegaard splittings by cutting along a Haken sphere. \\

\subsection*{Heegaard diagrams}
A genus $g$ Heegaard splitting can be represented as a standard closed orientable genus $g$ surface in $S^3$ with a pair of $g$ essential curves on it as follows: 
Identify the Heegaard surface $\Sigma$ with a standard closed orientable surface $S$ in $S^3$. 
Let $D_1$, \dots , $D_g$ be disks in one handlebody separated by $\Sigma$ such that they cut this handlebody into a 3-ball and $E_1$, \dots , $E_g$ be disks 
in the other handlebody such that they cut it into a 3-ball. 
Then we get a pair of $g$ curves, {$\partial D_1$, \dots , $\partial D_g$} and {$\partial E_1$, \dots , $\partial E_g$}. 
Draw this curves on $S$. This presentation is called a {\it Heegaard diagram}. 
We can reconstruct a Heegaard splitting by a Heegaard diagram: Paste disks along one of pair of curves in the interior of the standard surface, 
paste disks along the other curves in the exterior of the standard surface and paste two 3-balls along remaining spheres. \\

\subsection*{Fibered links}
Let $M$ be a closed orientable 3-manifold and $L$ be a link in $M$. $L$ is called a fibered link if $Cl(M \setminus N(L))$ is a fiber bundle over $S^{1}$ 
whose fiber is an orientable surface and the boundary of each fiber is isotopic to $L$ in $N(L)$, 
where $N(L)$ is a regular closed neighborhood of $L$ in $M$ and $Cl(\cdot )$ is the closure. 
If $L$ is a knot and the fiber is a torus $T$ with one boundary component, $L$ is called a GOF-knot in $M$. Let $L$ be a fibered link in $M$ and $F$ be its fiber. 
By thickening $F$, $Cl(M \setminus N(L))$ is decomposed into two handlebodies of the same genus, $g$. 
Moreover, by the property of fibered links we get a genus $g$ Heegaard splitting of $M$ 
such that $L$ is on its Heegaard surface and $L$ decomposes each handlebody to $F\times I$.
In particular, if $K$ is a GOF-knot in $M$, there is a genus two Heegaard splitting of $M$ and $K$ is on the Heegaard surface, 
decomposing each handlebody to $T\times I$. 
In this paper, two fibered links $L_1$ and $L_2$ in $M$ are said to be equivalent if there is a fiber preserving self-homeomorphism of $M$ sending $L_1$ to $L_2$.\\

\subsection*{Plumbing}
Let $L_1$ and $L_2$ be two fibered links in two closed orientable 3-manifolds $M_1$ and $M_2$ respectively. 
Let $F_i$ be a fiber of $L_i$ in $M_i$ ($i=0,1$). Then we can construct a fibered link in $M_1\# M_2$ from $L_1$ and $L_2$ as follows: 
Let $\alpha _i$ be a properly embedded essential arc in $F_i$ and $D_i$ be a small closed neighborhood of $\alpha _i$ in $F_i$. 
$D_i$ can be identified with $\alpha _{i}\times [-1,1]$. 
We construct a new surface $F$ by pasting $D_1$ and $D_2$ so that for every $t\in [-1,1]$, $\alpha _{1}\times \{t\}$ is identified with 
an arc intersecting once to $\alpha _2\times \{s\}$ for every $s\in [-1,1]$. 
For such an operation, we say that $F$ is obtained by the $plumbing$ of $F_1$ and $F_2$. 
In \cite{14}, Stallings showed that a surface obtained by the plumbing of two surfaces which are 
fibers of two fibered links in $S^3$ is also a fiber of a fibered link in $S^3$. 
This statement can be generalized to arbitrary closed orientable 3-maniflds. Thus, $\partial F$ is a fibered link in $M_1\# M_2$ with $F$ as a fiber surface.\\

\subsection*{Monodromy}
Let $L$ be a fibered link in a closed orientable 3-manifold $M$ and $F$ be a fiber surface of $L$. Then $Cl(M \setminus N(L))$ is a $F$-bundle over $S^1$. 
It is obtained from $F\times [0,1]$ by pasting $F\times \{0\}$ and $F\times \{1\}$ using an orientation preserving self-homeomorphism of $F$. 
This map is called the monodromy of $L$ (and $F$) or the monodromy of $Cl(M \setminus N(L))$. 
Let $f$ and $g$ be two orientation preserving self-homeomorphisms of $F$. 
Note that under orientation- and fiber preserving homeomorphisms, $F$-bundles over $S^1$ using $f$ and $g$ are equivalent if and only if 
$f$ and $g$ are in the same conjugacy class of the mapping class group of $F$ (each component of $\partial F$ is fixed setwise). 
If we work under fiber preserving homeomorphisms, an orientation reversing map can be added. 
In particular the monodromy of a GOF-knot is classified in the conjugacy classes in $GL_{2}(\mathbb{Z})$. 
In fact we can say about monodromies under the plumbing as followings: 
Let $F_1$ and $F_2$ are fibers of two fibered links in two closed orientable 3-manifolds $M_1$ and $M_2$ respectively, 
and let $f_1$ and $f_2$ be monodromies of $F_1$ and $F_2$. 
Then the monodromy of $F$, obtained by the plumbing of $F_1$ and $F_2$ is $\tilde{f_{1}}\circ \tilde{f_{2}}$ where $\tilde{f_i}$ is an extension of $f_i$ to $F$. \\

\subsection*{Manifolds which have genus one Heegaard splittings}
Every Heegaard diagram of a genus one Heegaard splitting is a standard torus with two simple closed curves on it, 
one is the meridian curve and the other is a $(p,q)$-curve ($p$ and $q$ are coprime). 
We may assume $p$ is non-negative. If $(p,q)=(1,0)$ or $(1,\pm 1)$, the manifold is $S^3$. If $(p,q)=(0,\pm 1)$, the manifold is $S^2\times S^1$. 
If otherwise, the manifold is called the lens space of type $(p,q)$, $L(p,q)$. 
For $(p,q)$, we set $q'$ to be the least non-negative number such that $qq'\equiv 1$ mod $p$. 
Note that $L(p,q)$ is homeomorphic to $L(p,q')$ (by changing the roles of two handlebodies of the genus one Heegaard diagram). 
It is known that the genus one Heegaard splitting of $S^3$, $S^2\times S^1$ or $L(p,q)$ is unique under homeomorphisms \cite{4}, \cite{15}.\\

\subsection*{Fibered annulus}
A fibered annulus is an annulus in a closed orientable 3-manifold $M$ which is a fiber of a fibered link. 
As above, a manifold which has a fibered annulus has a genus one Heegaard splitting. 
Moreover because of the fact that the group of self-homeomorphisms of an annulus is generated by the Dehn twist along the core curve 
and the properties of fibered links, the corresponding genus one Heegaard diagram is a standard torus with the meridian curve and ($p,\pm 1$)-curve. 
Hence a manifold which has a fibered annulus is $S^3$, $S^2\times S^1$ or $L(p,\pm 1)$. 
We can see easily that each of $S^2\times S^1$, $L(p,1)$ and $L(p,-1)$ ($p\neq 2$) has just one fibered annulus 
and each of $S^3$, $L(2,1)$ and $L(2,-1)$ has just two fibered annuli under orientation preserving self-homeomorphisms. 
(Note that there is an orientation reversing homeomorphism between $L(2,1)$ and $L(2,-1)$ as noted in \cite{10}.)
The monodromy of the fibered annulus in $S^2\times S^1$ is the identity map, 
that in $L(p,1)$ (so called \emph{$p$-Hopf band} in \cite{2}) is the $p$ times positive Dehn twists along the core curve ($p\neq 2$), 
that in $L(p,-1)$ (so called \emph{$-p$-Hopf band} in \cite{2}) is the $p$ times negative Dehn twists along the core curve ($p\neq 2$), 
that of one fibered annulus in $L(2,1)$ (so called \emph{$2$-Hopf band} in \cite{2}) is the 2 times positive Dehn twists along the core curve, 
that of the other fibered annulus in $L(2,1)$ (so called \emph{$-2$-Hopf band} in \cite{2}) is the 2 times negative Dehn twists along the core curve, 
that of the $+1$-Hopf annulus (resp. $-1$-Hopf annulus) in $S^3$ is the positive (resp. negative) Dehn twist along the core curve.\\

\subsection*{A genus two Heegaard splitting of a manifold which admits a reducible one.}
If a closed orientable 3-manifold $M$ admits a reducible genus two Heegaard splitting, then $M$ is homeomorphic to 
$(S^2\times S^1)\# (S^2\times S^1)$, $(S^2\times S^1)\# L(p,q)$, $L(p_1,q_1)\# L(p_2,q_2)$, $S^2\times S^1$, $S^3$, or $L(p,q)$. 
In \cite{6}, Casson and Gordon proved that if a 3-manifold is reducible (i.e. has a sphere which does not bound a 3-ball), every Heegaard splitting of it is reducible. 
In \cite{15}, Waldhausen proved that every Heegaard splitting of $S^3$ whose genus is greater than $0$ is reducible, 
and in \cite{4}, Bonahon-Otal proved that every Heegaard splitting of lens spaces whose genus is greater than $1$ is reducible. 
These imply that every genus two Heegaard splitting of the above manifolds is reducible. 
As mentioned above, every reducible genus two Heegaard splitting can be decomposed into two genus one Heegaard splittings. 
In the opposite direction, every reducible genus two Heegaard splitting is obtained by connecting two genus one Heegaard splittings. 
Hence the manifold which has a reducible genus two Heegaard splitting has at most two genus two Heegaard splittings under homeomorphisms. 
It depends on the choice of a solid torus of one genus one Heegaard splitting which should be connected to one solid torus of the other genus one Heegaard splitting.  
In \cite{8}, Cho and Koda gave the necessary and sufficient condition for a Heegaard splitting of such manifolds of being unique under homeomorphisms. 
(Unique in $S^3$, $L(p,q)$, $S^2\times S^1$, $(S^2\times S^1)\# (S^2\times S^1)$, $L(p,q)\# (S^2\times S^1)$ and 
$L(p_1,q_1)\# L(p_2,q_2)$ with ${q_1}^2 \equiv 1$ mod $p_1$ or ${q_2}^2 \equiv 1$ mod $p_2$ and not unique in the other part of $L(p_1,q_1)\# L(p_2,q_2)$.) 
However, even though there are two, their Heegaard diagrams are very similar. 
We will focus on one Heegaard surface and its diagram. The arguments for the other Heegaard splittings are similar and the conclusions are the same.\\\vspace{0.5in}

\ From the above, we see that if $M$ has a GOF-knot, then $M$ must have a genus two Heegaard splitting and that if $M$ has a reducible genus two Heegaard splitting 
then $M$ must be homeomorphic to $(S^2\times S^1)\# (S^2\times S^1)$, $(S^2\times S^1)\# L(p,q)$, $L(p_1,q_1)\# L(p_2,q_2)$, $S^2\times S^1$, $S^3$, or $L(p,q)$. 
(Moreover, the genus two Heegaard splitting obtained by the GOF-knot is also reducible.) 
Therefore $M$ which admits a reducible genus two Heegaard splitting and possibly has a GOF-knot is 
$(S^2\times S^1)\# (S^2\times S^1)$, $(S^2\times S^1)\# L(p,q)$, $L(p_1,q_1)\# L(p_2,q_2)$, $S^2\times S^1$, $S^3$, or $L(p,q)$. \\\\

\ In this paper, our goal is to give another proof of the following theorems already obtained in \cite{2}, \cite{3}, \cite{12}.\\
 \begin{thm} \label{thm1}
  $\rm{(Baker \cite{2})}$ \ There is just one GOF-knot in $(S^2 \times S^1) \# (S^2 \times S^1)$
  and the fiber is obtained by the plumbing of two fibered annuli in $S^2\times S^1$.
 \end{thm}

 \begin{thm} \label{thm2}
  $\rm{(Baker \cite{2})}$ \ There is a GOF-knot in $L(p,q)\# (S^2\times S^1)$ if and only if $q \equiv \pm 1$ mod $p$,
  and if there is a GOF-knot, the fiber can be obtained by the plumbing of each of fibered annulus in $L(p,q)$ and $S^2\times S^1$. 
 \end{thm}

 \begin{thm} \label{thm3}
  $\rm{(Baker \cite{2})}$ \ There is a GOF-knot in $L(p_1,q_1)\# L(p_2,q_2)$ if and only if $q_i \equiv \pm 1$  mod $p_i$ ($i=1,2$),
  and if there is a GOF-knot, the fiber can be obtained by the plumbing of fibered annuli in $L(p_1,q_1)$ and $L(p_2,q_2)$. 
 \end{thm}

 \begin{thm} \label{thm4}
  $\rm{(Morimoto \cite{12})}$ \ There is just one GOF-knot in $S^2\times S^1$.
 \end{thm}

 \begin{thm} \label{thm5}
  $\rm{(Morimoto \cite{12})}$ \ There are just two GOF-knots in $S^3$. They are the $($left and right hand $)$ trefoil and the figure eight knot.
 \end{thm}

 \begin{thm} \label{thm6}
  $\rm{(Baker \cite{3})}$ \\
  $\rm{(1)}$ $L(4,1)$ has just three GOF-knots.\\
  $\rm{(2)}$ $L(p,1)$ ($p \neq 4$) has just two GOF-knots.\\
  $\rm{(3)}$ For positive integers $a$ and $b$, $L(2ab+a+b+1,2b+1)$ has just one GOF-knot.\\
  $\rm{(4)}$ For positive integers $a$ and $b$, $L(2ab+a+b,2b+1)$ (except for $L(4,3)$) has just one GOF-knot.\\
  $\rm{(5)}$ A lens space which is not homeomorphic to any of the above types has no GOF-knots.
 \end{thm}
\vspace{0.3in}

\ Morimoto proved some of the above theorems by using the monodromy of a GOF-knot in a 3-manifold and the fundamental group of this manifold. 
Baker proved the other theorems by using the one-to-one correspondence between GOF-knots and axes of closed 3-braids in $S^3$. 
(There is a double branched covering map branched along a closed 3-braid. The preimage of the axis is a GOF-knot.)
In this paper, we give a unified proof of the above theorems by using the fact that the genus two Heegaard splitting of these manifolds is almost unique. 
As a result, we find the positions of GOF-knots clearly on standard Heegaard surfaces of these manifolds.
\vspace{0.8in}

\section{Methods} \label{sec3}
\begin{figure}[htbp]
 \begin{center}
  \includegraphics[width=50mm]{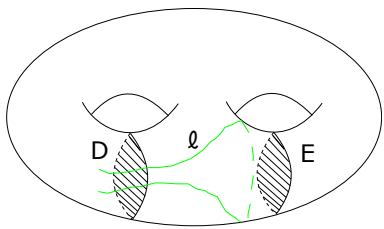}
 \end{center}
 \caption{uncancellable $xx^{-1}$}
 \label{fig:one}
\end{figure}

\begin{figure}[htbp]
 \begin{center}
  \includegraphics[width=40mm]{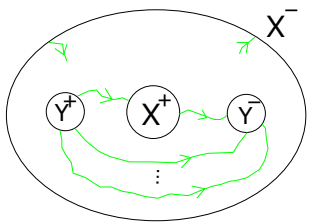}
 \end{center}
 \caption{$xy^{n}x^{-1}$}
 \label{fig:two}
\end{figure}

\begin{figure}[htbp]
 \begin{center}
  \includegraphics[width=60mm]{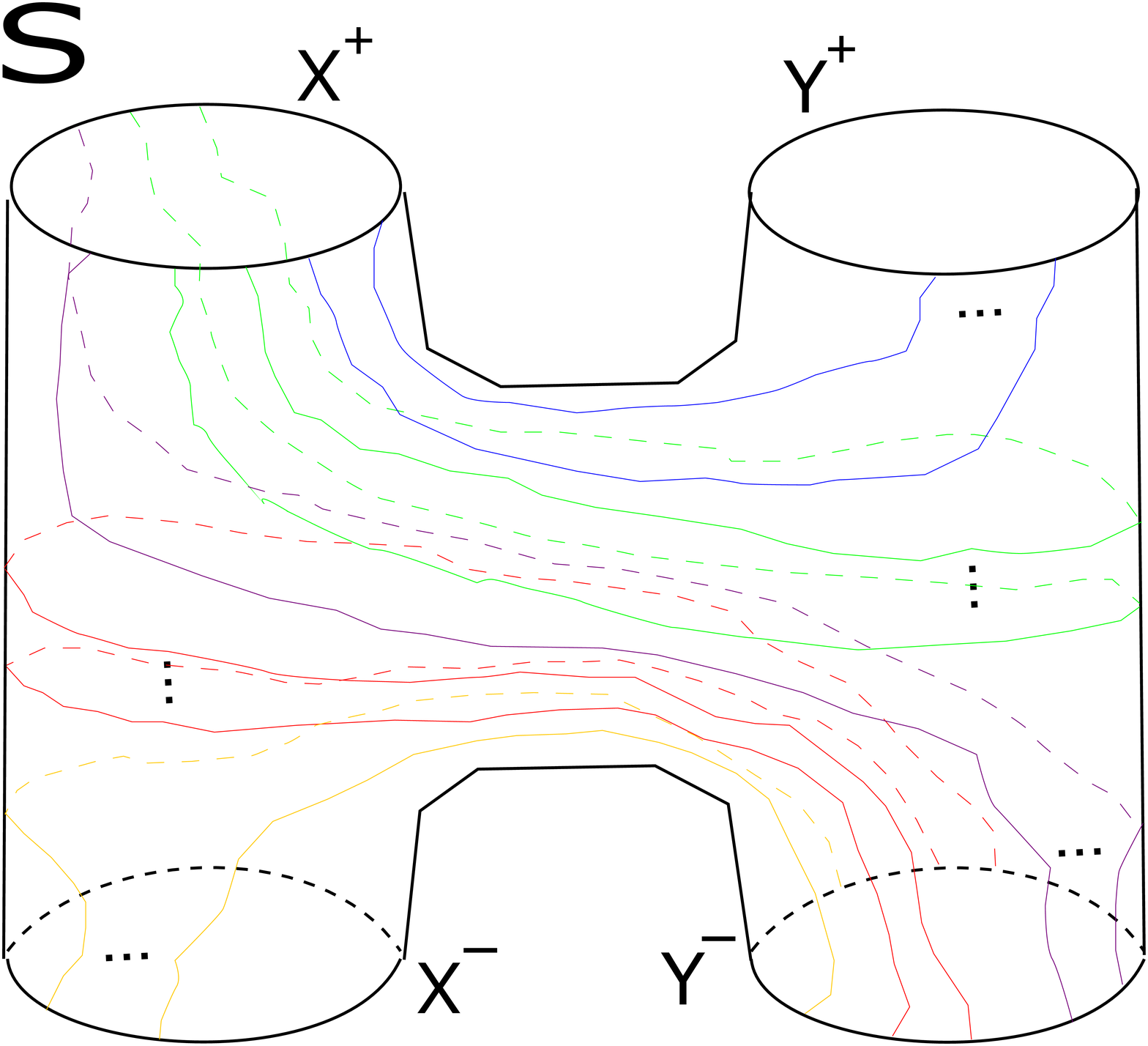}
 \end{center}
 \caption{possible 4 types of arcs when there are $xx^{-1}$ and $y^{-1}y$}
 \label{fig:three}
\end{figure}

\begin{figure}[htbp]
 \begin{center}
  \includegraphics[width=60mm]{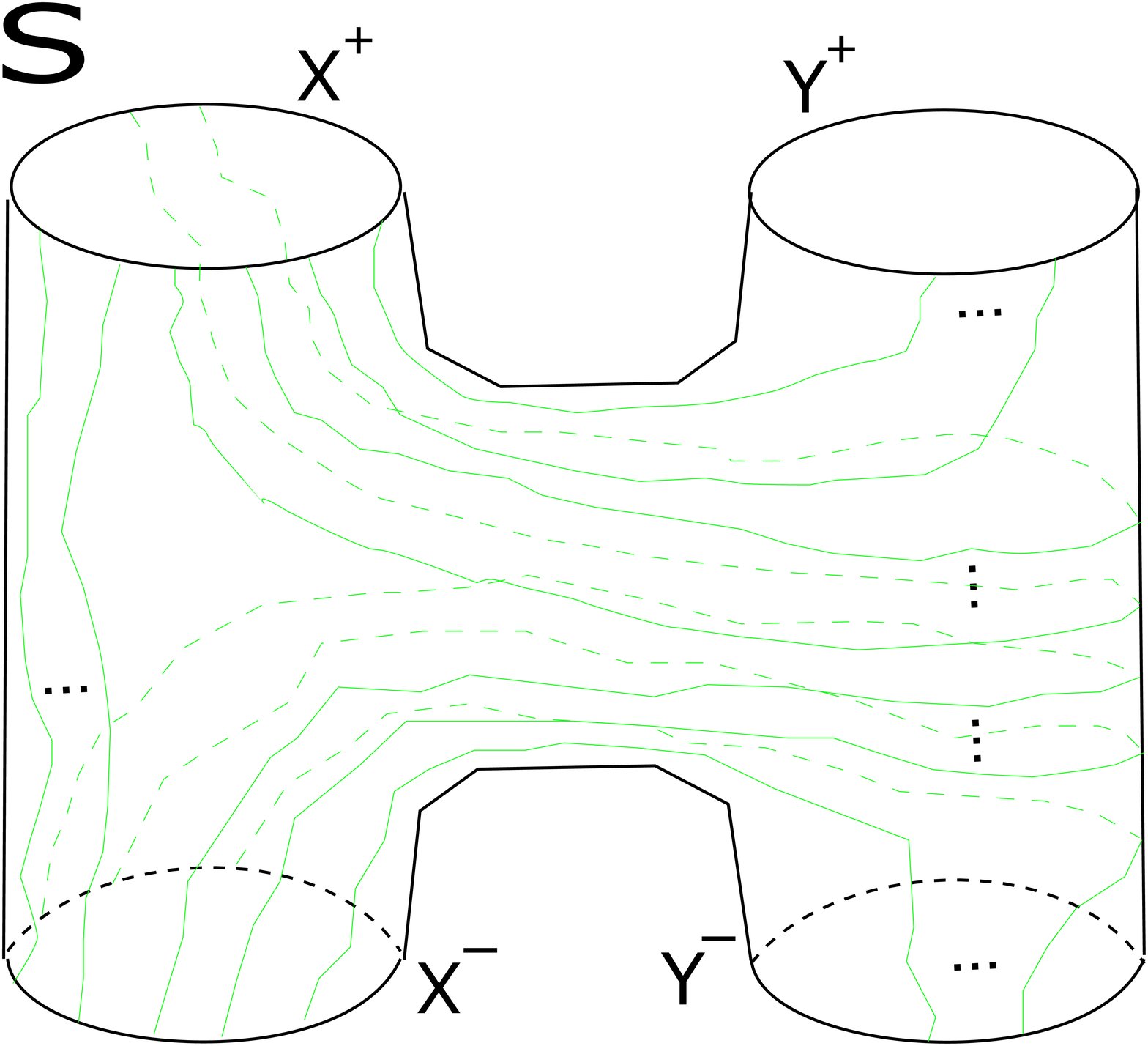}
 \end{center}
 \caption{actually possible arcs when there is $xx^{-1}$}
 \label{fig:four}
\end{figure}

\begin{figure}[htbp]
 \begin{minipage}{0.5\hsize}
  \begin{center}
   \includegraphics[width=60mm]{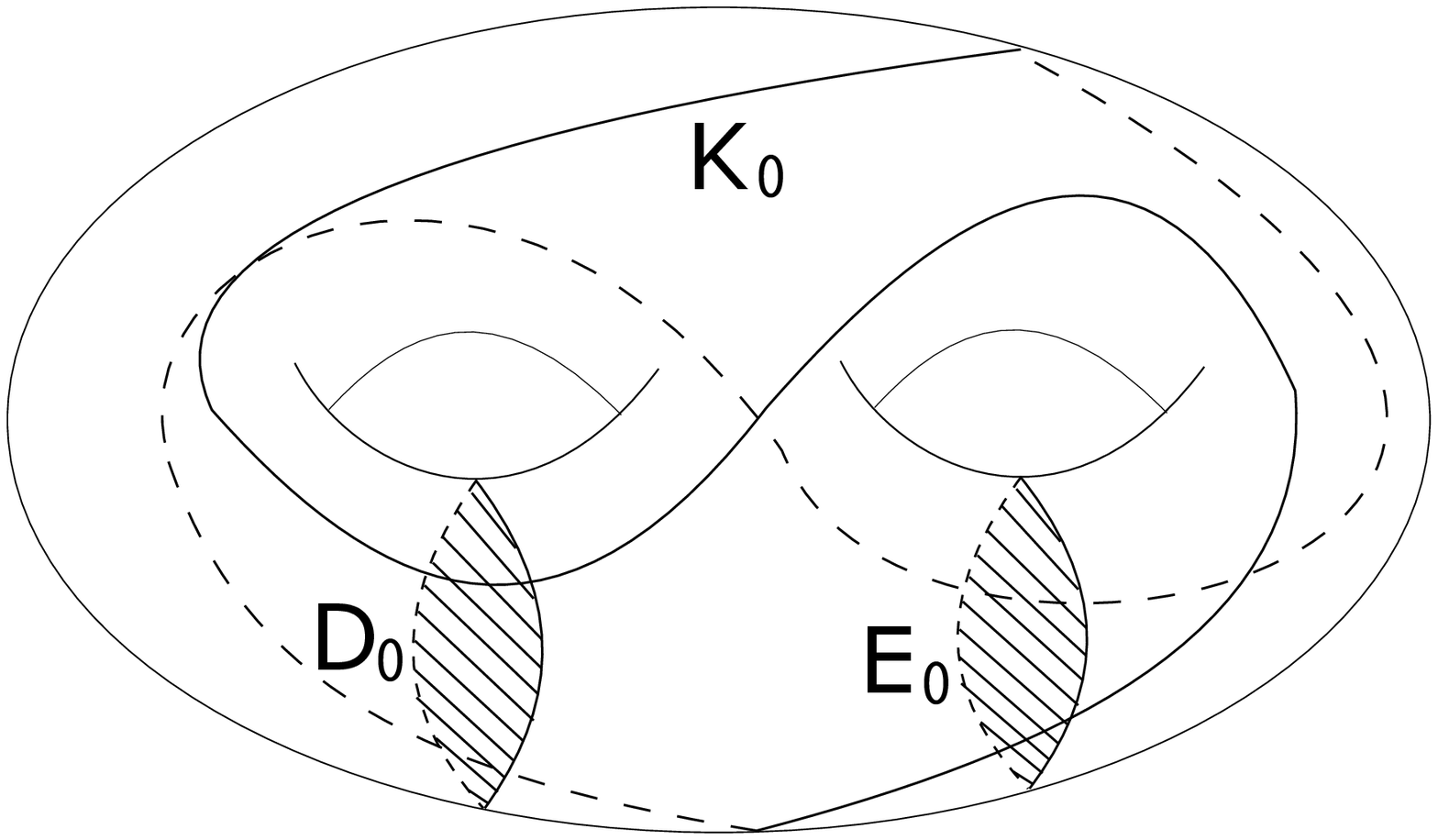}
  \end{center}
  \caption{$\partial D_0$, $\partial E_0$ and $K_0$}
  \label{fig:five}
 \end{minipage}
 \begin{minipage}{0.5\hsize}
  \begin{center}
   \includegraphics[width=45mm]{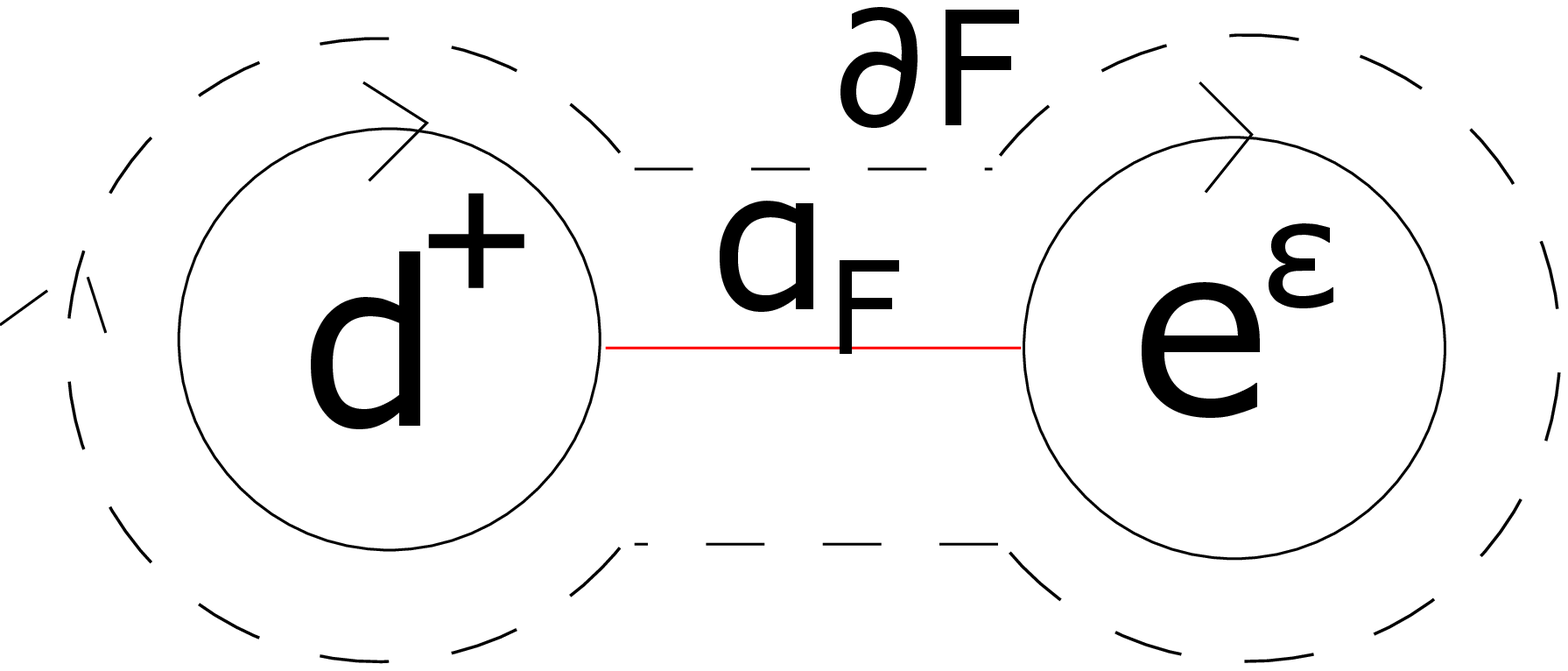}
  \end{center}
  \caption{$\alpha_F$ and $F$}
  \label{fig:six}
 \end{minipage}
\end{figure}

\begin{figure}[h]
 \begin{minipage}{0.5\hsize}
  \begin{center}
   \includegraphics[width=50mm]{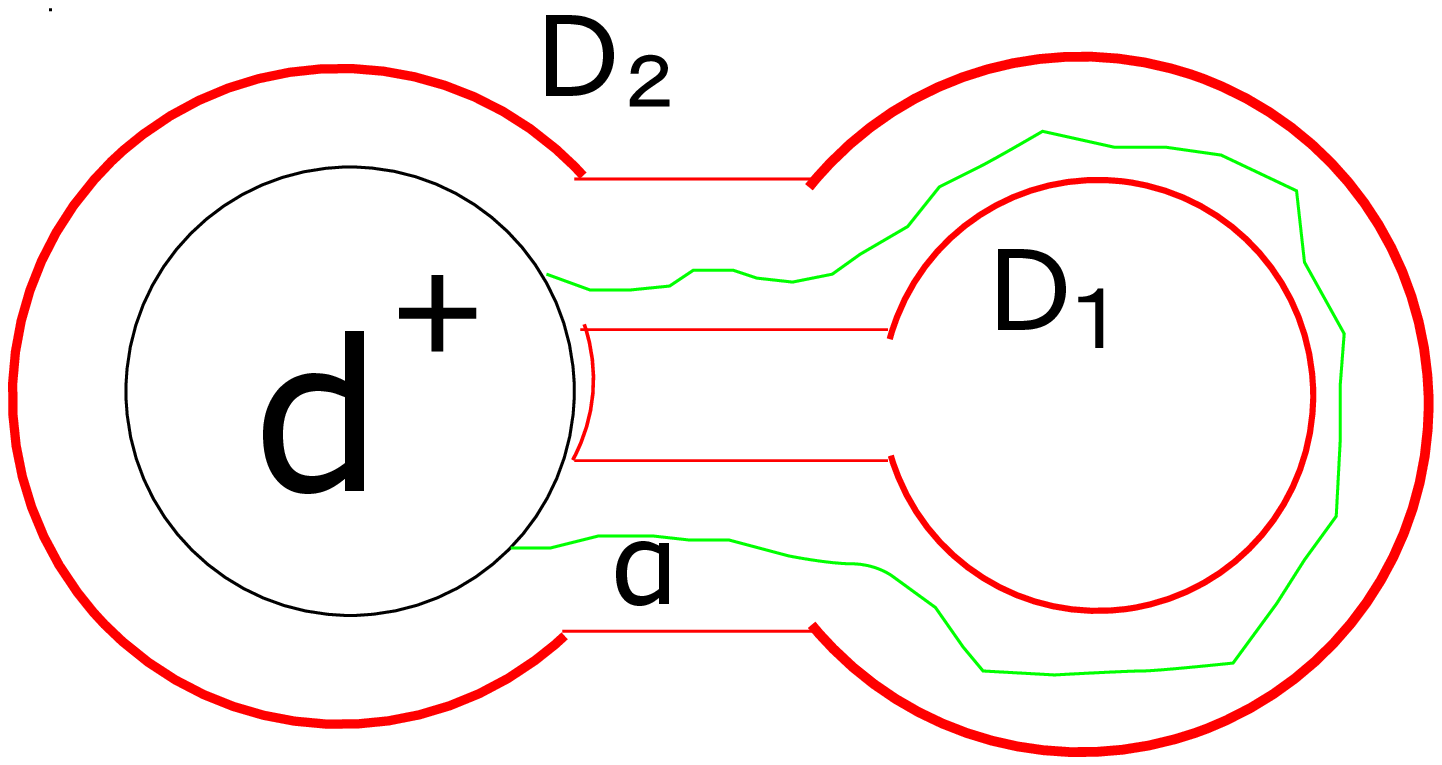}
  \end{center}
  \caption{$D$, $D_1$ and $D_2$}
  \label{fig:seven}
 \end{minipage}
 \begin{minipage}{0.5\hsize}
  \begin{center}
   \includegraphics[width=60mm]{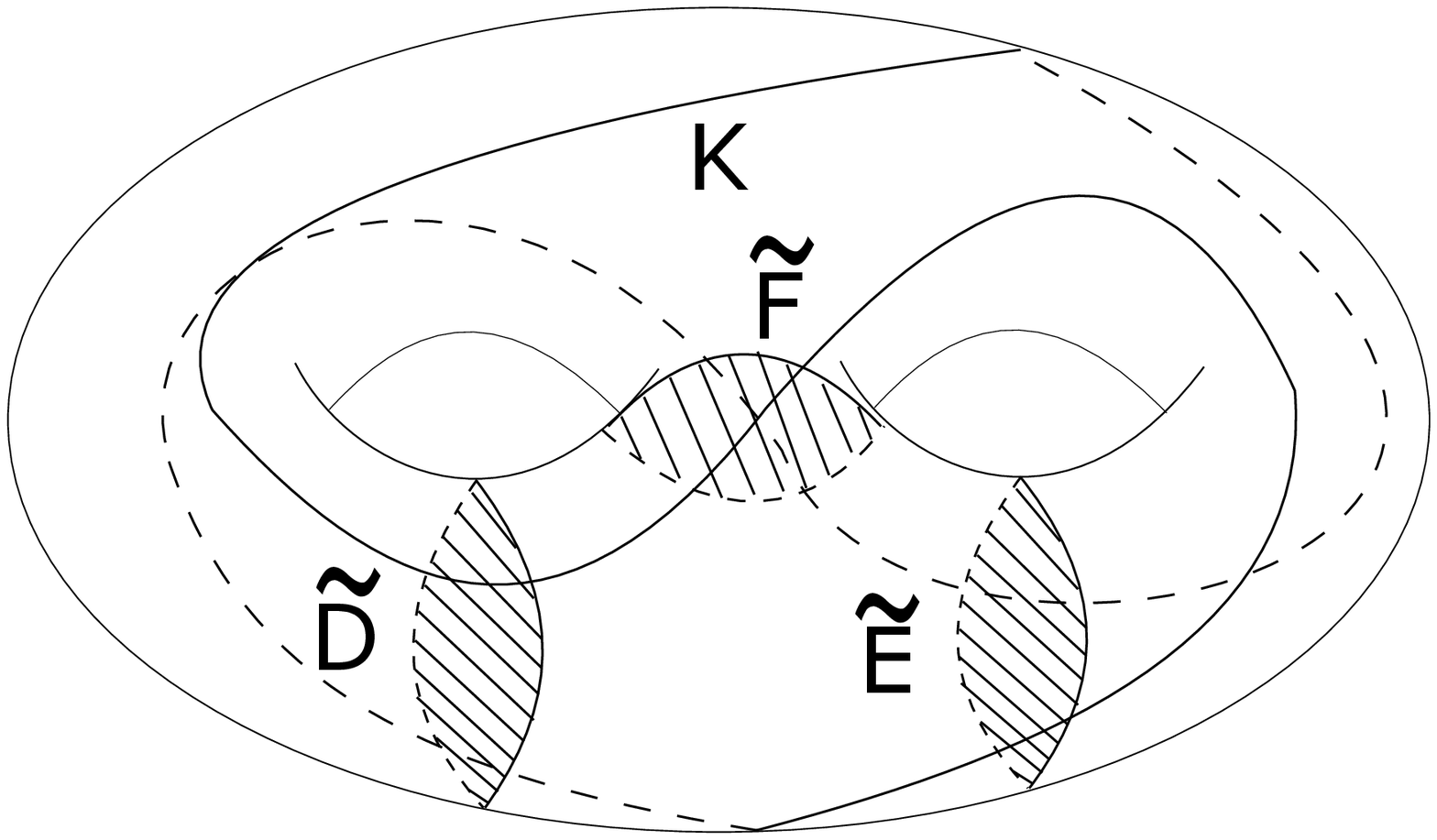}
  \end{center}
  \caption{$\partial \tilde{D}$, $\partial \tilde{E}$, $\partial \tilde{F}$ and $K$}
  \label{fig:eight}
 \end{minipage}
\end{figure}

As in Section \ref{sec2}, every GOF-knot of $M$ is on a genus two Heegaard surface so that it decomposes each handlebody into $T\times I$ ($T$ is a fiber surface). 
We call a simple closed curve on the boundary of genus two handlebody such that it decomposes the handlebody into $T\times I$ a \emph{GOF-knot on the handlebody}. 
To investigate the properties of such curves, we consider simple closed curves on the boundary of a genus two handlebody. \\\\

\ Let $V$ be a genus two handlebody and let $D$ and $E$ be two disjoint properly embedded, non-separating disks in $V$ which are not parallel. 
Then $D$ and $E$ cut $V$ into a 3-ball. 
Fix orientations of $\partial D$ and $\partial E$ and give them letters $x$ and $y$, respectively.
Let $l$ be an oriented simple closed curve on $\partial V$. Isotope $l$ so that $l$ intersects $\partial D \cup \partial E$ minimally and transversely. 
Then $l$ determines a word of $x$ and $y$ that can be read off by the intersections of $l$ with $\partial D$ and $\partial E$ in taking orientations into account. 
This word is well-defined up to cyclic permutations. Note that this word may not be reduced (i.e. may have subword of type $xx^{-1}$. See Figure \ref{fig:one}.). 
For a simple closed curve $l$, the word should be written cyclically, however we write it simply not cyclically. 
For this reason, we say that the word is reduced if it is cyclically reduced. The following lemma is frequently used later. \\
\\
\begin{lem} \label{lem1}
 $\rm{(Cho, Koda, \cite{7} \ \cite{8} \ \cite{9})}$\ In the above setting, if the word determined by a simple closed curve $l$ contains
 a subword of the form $xy^{n}x^{-1}$ for some $n \in \mathbb{N}$, then this word is reduced.
\end{lem}
{\it Proof from \rm{\cite{7}}}. \ Let $S$ be a sphere with four boundary components which is obtained from $\partial V$ by cutting along $\partial D \cup \partial E$. 
We denote by $X^{+}$, $X^{-}$, $Y^{+}$ and $Y^{-}$ the boundary components of $S$ coming from $\partial D$ and $\partial E$, respectively. See Figure \ref{fig:two}. 
In $S$, the subword $xy^{n}x^{-1}$ corresponds to $n+1$ arcs. 
The first one starts from $X^{+}$ and ends at $Y^{-}$, the second one starts from the point of $Y^{+}$ 
corresponding to the terminal point of the first one in $\partial V$ and ending at $Y^{-}$, ..., the $n$-th one starts from $Y^{+}$ and ends at $Y^{-}$, 
the last one starts from $Y^{+}$ and ends at $X^{+}$. 
Then there must be two arcs so that one connects $X^{-}$ and $Y^{+}$, the other connects $X^{-}$ and $Y^{-}$. 
(For example, the next arc of the above subarcs of $l$ starts at $X^{-}$ and ends at $X^{+}$, $Y^{+}$ or $Y^{-}$ because of minimality. 
If the ending point is on $X^{+}$, the next arc starts from $X^{-}$ and ends at $X^{+}$, $Y^{+}$ or $Y^{-}$. 
Repeating this in finitely many times, we get an arc connecting $X^{-}$ with $Y^{+}$ or $Y^{-}$.) 
It follows that we cannot draw an arc of the form $xx^{-1}$ or $yy^{-1}$ without intersecting $l$. 
Therefore $l$ cannot contain a subword of the form $xx^{-1}$ or $yy^{-1}$. This implies the word represented by $l$ is reduced. 
\qed
\\
\ The followings are implicitly used later. 
\begin{lem} \label{lem2}
 \ In the above setting, if a simple closed curve $l$ on $\partial V$ has a subarc representing $xx^{-1}$ $($or $x^{-1}x$$)$,
 then $l$ has no subarcs representing $yy^{-1}$ nor $y^{-1}y$.
\end{lem}
{\it Proof}. \ Let $S$ be a sphere with four boundary components which is obtained from $\partial V$ cutting along $\partial D \cup \partial E$. 
We denote by $X^{+}$, $X^{-}$, $Y^{+}$ and $Y^{-}$ the boundary components of $S$ coming from $\partial D$ and $\partial E$, respectively. 
For a subarc $c$ of $l$ representing $xx^{-1}$, there is an arc $\alpha$ in $S$ connecting $X^{+}$ and another boundary component $W$ such that 
$c$ is the boundary of a regular neighborhood of $\alpha \cup W$ in the interior of $S$, $Int(S)$. 
If $W$ was $X^{-}$, the number of the points of $l \cap X^{-}$ does not coincide with that of $l \cap X^{+}$. It cannot occur. 
Thus $W$ is $Y^{+}$ or $Y^{-}$. We assume $W$ is $Y^{+}$. In this case, there are no subarcs of $l$ representing $yy^{-1}$. 
If there is a subarc of $l$ representing $y^{-1}y$, the subarc $c$ and this arc on $S$ are of the form in Figure \ref{fig:three}. 
Let $a$, $b$, $c$, $d$ and $e$ be the number of subarcs of $l$ on $S$ connecting $X^{+}$ and $X^{+}$, 
connecting $X^{+}$ and $Y^{+}$, connecting $X^{+}$ and $Y^{-}$, connecting $Y^{-}$ and $Y^{-}$ and connecting $Y^{-}$ and $X^{-}$. 
Note that there are no subarcs of the other types (See Figure \ref{fig:three}). 
$a$ and $d$ are at least $1$. Then two equalities $2a+b+c=e$ and $b=c+2d+e$ must hold. It cannot occur. Therefore there are no subarcs of $l$ representing $y^{-1}y$. 
\qed
\\
\begin{lem} \label{rem}
 \ In the same setting of Lemma \ref{lem2}, if there are $n$ subarcs of simple closed curve $l$ on $\partial V$ each of which represents $xx^{-1}$
 then there must be $n$ subarcs of simple closed curve $l$ on $\partial V$ each of which represents $x^{-1}x$.
\end{lem}
{\it Proof}. It follows by an argument similar to Lemma \ref{lem2}.
 $l$ is like in Figure \ref{fig:four}. 
\qed
\\
\\\\
\begin{lem} \label{lem3}
 \  Let $D$ and $E$ be two disjoint properly embedded, non-separating disks which are not parallel and $K$ be an oriented GOF-knot on a genus two handlebody $V$.
 Then by assigning $x$ and $y$ to $\partial D$ and $\partial E$ with an appropriate orientation,
 the word represented by $K$ becomes the commutator $xyx^{-1}y^{-1}$ of $x$ and $y$ after reduction.
\end{lem}
{\it Proof}. At first, consider the GOF-knot ${}K_0$ and the disjoint properly embedded, non-separating disks $D_0$ and $E_0$ 
which are not parallel in Figure \ref{fig:five}. It is easy to orient ${}\partial D_0$ and ${}\partial E_0$ with an assignment letters $x$ and $y$ to them 
so that $K_0$ represents the commutator $xyx^{-1}y^{-1}$.\\
\ Let $D$, $E$ and $F$ be disjoint properly embedded, non-separating disks which are not pairwise parallel. 
We will show that if the word represented by ${}K_0$ becomes the commutator $xyx^{-1}y^{-1}$ after reduction 
by giving $\partial D$ and $\partial E$ an appropriate orientation and letters $x$ and $y$, 
then the word becomes the commutator $zwz^{-1}w^{-1}$ after reduction by giving $\partial D$ and $\partial F$ an appropriate orientation and letters $z$ and $w$. 
Let $\Sigma$ denote $\partial V$, ${}\Sigma'$ denote $\Sigma \setminus \partial D \cup \partial E$, and $d^{+}$, $d^{-}$, $e^{+}$, $e^{-}$ 
denote the boundary components of ${}\Sigma '$ coming from $\partial D$ and $\partial E$. 
There is an arc ${}\alpha _{F}$ in ${}\Sigma '$ connecting ${}d^{+}$ and ${}e^{\epsilon}$ ($\epsilon \in \{\pm \}$), 
such that one boundary component of a small neighborhood of $d^{+} \cup {}\alpha_{F} \cup e^{\epsilon}$ is (isotopic to) $\partial F$. 
(the others are $d^{+}$ and $e^{\epsilon}$.)
Tentatively we give $\partial F$ the orientation coming from $\partial D$ (see Figure \ref{fig:six}). 
Isotope ${}\alpha _F$ so that its endpoints are not on $K_0$. 
The word represented by the subarc of $K_0$ cut by $\partial F$ near every intersection point of $K_0$ with $\alpha _F$ is $ww^{-1}$ or $w^{-1}w$. 
Thus up to reduction, the word represented by $K_0$ in letters $z$ and $w$ comes from the intersections of $\partial D$ and $\partial E$ with ${}K_0$. 
The small subarc of $K_0$ representing the word $x$ in letters $x$ and $y$ corresponds to the word $zw$ in letters $z$ and $w$. 
Similarly, $x^{-1}$ corresponds to $w^{-1}z^{-1}$, $y$ corresponds to $w$ or $w^{-1}$ 
(depending on the choice of the orientation of $\partial F$) and $y^{-1}$ corresponds to $w^{-1}$ or $w$. 
Therefore if the word represented by ${}K_0$ in $x$ and $y$ is $xyx^{-1}y^{-1}$ after reduction, 
the word in $z$ and $w$ is $zwz^{-1}w^{-1}$ or $zw^{-1}z^{-1}w$ after reduction. 
Changing the orientation of $\partial F$ if necessary, the word is $zwz^{-1}w^{-1}$.\\
\\
\ Because of the connectivity of the non-separating disk complex of a handlebody \cite{11}, 
every pair of disjoint properly embedded non-separating disks $D$ and $E$ in $V$ 
which are not parallel can be constructed by an iteration of the above operation to $D_0$ and $E_0$. 
Therefore for any disjoint two properly embedded, non-separating disks $D$ and $E$ in $V$ which are not parallel, 
the word represented by ${}K_0$ is $xyx^{-1}y^{-1}$ after reduction by giving the boundaries of them an appropriate orientation and letters $x$ and $y$.\\
\ For any GOF-knot $K$, there is a self-homeomorphism $f$ of $V$ such that $f$ sends $K$ to ${}K_0$. 
The word of $K$ by using any disjoint two properly embedded, non-separating disks $D$ and $E$ in $V$ 
which are not parallel is the same as the word of ${}K_{0}=f(K)$ by using disks $f(D)$ and $f(E)$. 
By the procedure discussed above, it is the commutator after reduction.
\qed
\\\\\\
\ Conversely the following holds.\\
\begin{lem} \label{lem4}
 \ Let {\it D} and {\it E} be disjoint properly embedded, non-separating disks in genus two handlebody $V$
 which are not parallel and $K$ be a simple closed curve on $\partial V$.
 If the word represented by $K$ is $xyx^{-1}y^{-1}$ after reduction by giving $\partial D$ and $\partial E$ an appropriate orientation and letters $x$ and $y$,
 then $K$ is a GOF-knot on $V$.
\end{lem}
{\it Proof}. If there is a subarc $\alpha$ of $K$ representing a word $xx^{-1}$, $\alpha$ and an arc $\beta$ on $D$ 
which has common endpoints with $\alpha$ bound a disk $D'$ in $V$ whose interior is disjoint from $D$ and $E$. 
$\alpha$ cuts $D$ into two disks $D_{1}$ and $D_{2}$. At least one of $D' \cup D_{1}$ and $D' \cup D_{2}$ is a non-separating disk. 
Denote this disk by $\bar{D}$ (See Figure \ref{fig:seven}) and isotope $\bar{D}$ so that $\partial \bar{D}$ intersects $K$ minimally. 
Then by giving $\partial \bar{D}$ and $\partial E$ an appropriate orientation and letters $z$ and $w$ the number of letters of the word represented by $K$ 
in terms of $z$ and $w$ is less than that of the word in terms of $x$ and $y$. 
Since $\bar{D}$ is disjoint from and not parallel to $D$ and $E$, by the argument of Lemma \ref{lem3} the word represented by $K$ in terms of $z$ and $w$ is 
also the commutator of $z$ and $w$ after reduction. 
Using this operation repeatedly, we find two disjoint properly embedded, non-separating disks $\tilde{D}$ and $\tilde{E}$ in $V$ which are not parallel 
such that the word represented by $K$ is reduced and a form of $\tilde{x} \tilde{y} \tilde{x}^{-1} \tilde{y}^{-1}$ 
by giving $\partial \tilde{D}$ and $\partial \tilde{E}$ an appropriate orientation and letters $\tilde{x}$ and $\tilde{y}$. 
Let ${}\Sigma'$ denote $\partial V \setminus \partial \tilde{D} \cup \partial \tilde{E}$, and $\tilde{d}^{+}$, $\tilde{d}^{-}$, $\tilde{e}^{+}$, $\tilde{e}^{-}$ 
denote the boundary components of ${}\Sigma '$ coming from $\partial \tilde{D}$ and $\partial \tilde{E}$. 
The curve $K$ corresponds to four arcs on ${}\Sigma'$. The first one starts from $\tilde{d}^{+}$ and ends at $\tilde{e}^{-}$, 
the second starts from the point of $\tilde{e}^{+}$ corresponding to the terminal point of the first one in $\partial V$ and ends at $\tilde{d}^{+}$, 
the third starts from $\tilde{d}^{-}$ and ends at $\tilde{e}^{+}$, and the fourth starts from $\tilde{e}^{-}$ and ends at $\tilde{d}^{-}$. 
One boundary of a small regular neighborhood of $\tilde{d}^{+}$, $\tilde{e}^{-}$ and the first arc in ${}\Sigma'$ bounds a non-separating disk in $V$ 
which is disjoint from and not parallel to $\tilde{D}$ and $\tilde{E}$. Denote this disk by $\tilde{F}$. 
The simple closed curve $K$ and three disks $\tilde{D}$, $\tilde{E}$ and $\tilde{F}$ in $V$ are drawn in Figure \ref{fig:eight}. 
Though $K$ may be the mirror image of it and may be twisted along disks $\tilde{D}$, $\tilde{E}$ or $\tilde{F}$, 
we assume $K$ is like in Figure \ref{fig:eight} by a self-homeomorphism of $V$. 
The triplet ($K$, $\tilde{D}$, $\tilde{E}$) corresponds to the triplet ($K_{0}$, $D_{0}$, $E_{0}$) in Lemma \ref{lem3}. Therefore $K$ is a GOF-knot on $V$.
\qed
\\\\\\
\ By using the above two lemmas and the relation between GOF-knots and genus two Heegaard splittings, 
we get a one-to-one correspondence between a GOF-knot (with its fiber) and a simple closed curve on a genus two Heegaard surface 
whose representing words in both genus two handlebodies are commutators after reduction. 
The words can be read off by a Heegaard diagram.\\\\\\

\section{Individual cases} \label{sec4}
\subsection{$(S^2 \times S^1) \# (S^2 \times S^1)$}
\begin{figure}[htbp]
 \begin{minipage}{0.5\hsize}
  \begin{center}
   \includegraphics[width=50mm]{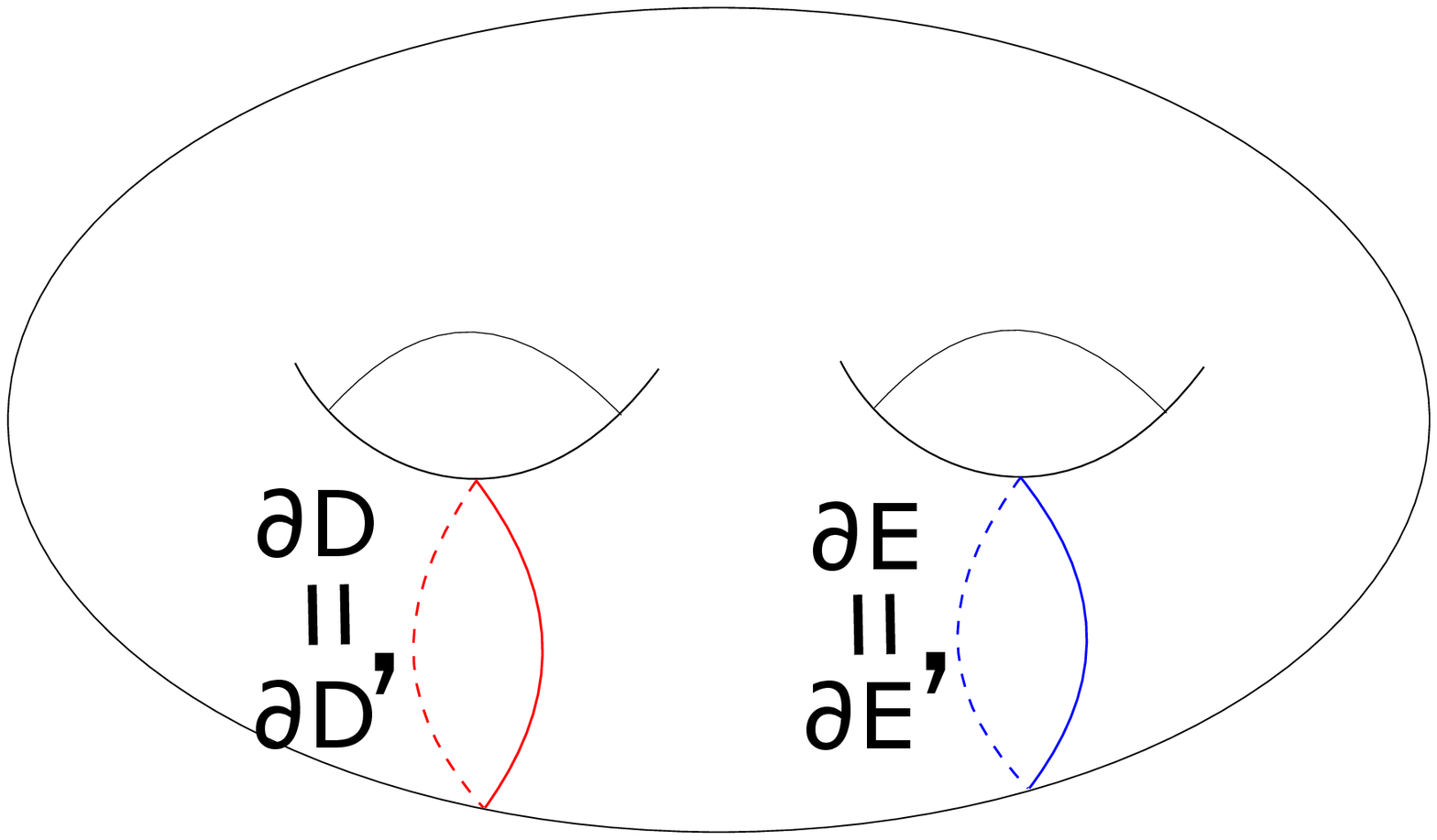}
  \end{center}
  \caption{a Heegaard diagram}
  \label{fig:nine}
 \end{minipage}
 \begin{minipage}{0.5\hsize}
  \begin{center}
   \includegraphics[width=50mm]{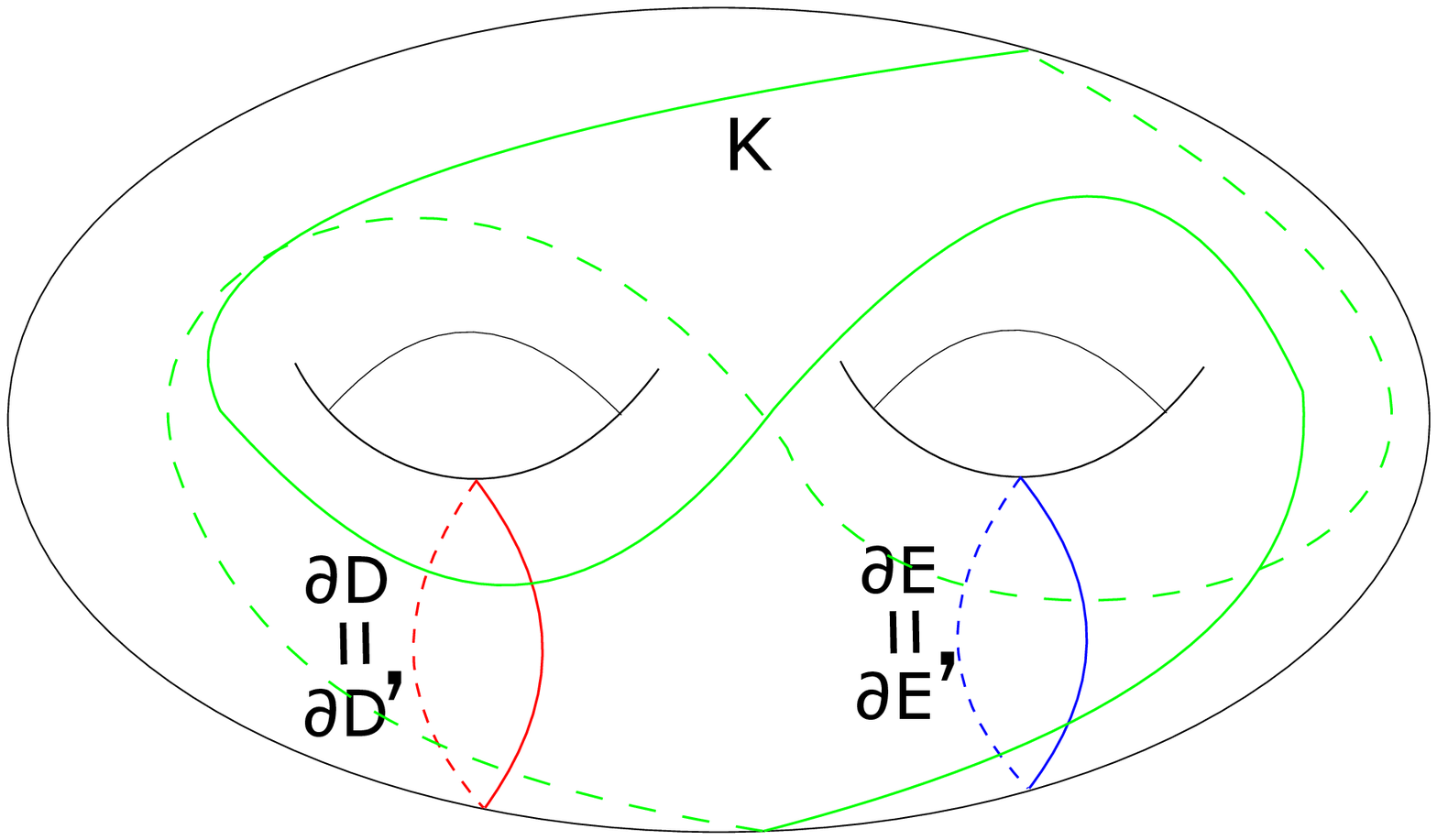}
  \end{center}
  \caption{one GOF-knot}
  \label{fig:ten}
 \end{minipage}
\end{figure}

In this case the genus two Heegaard splitting is unique as in Section \ref{sec2}. We consider a Heegaard diagram in Figure \ref{fig:nine}. 
We denote by $V \cup _{\Sigma} W$ the corresponding genus two Heegaard splitting. 
This case is exceptionally easy because of the following two properties.: 
One is that a GOF-knot in $V$ is also a GOF-knot in $W$ since $\partial D$ (resp. $\partial E$) is the same as $\partial D'$ (resp. $\partial E'$) in $\Sigma$, 
and the other is that the restriction on $\partial V = \Sigma$ of every self-homeomorphism of $V$ extends to a self-homeomorphism of $W$. 
By the first property, every GOF-knot in $(S^2 \times S^1) \# (S^2 \times S^1)$ corresponds to a GOF-knot in $V$. 
Since for any pair of two GOF-knots $K_1$ and $K_2$ in $V$, there is a fiber preserving self-homeomorphism of $V$ which takes $K_1$ to $K_2$. 
By the second property, all GOF-knots in $(S^2 \times S^1) \# (S^2 \times S^1)$ are equivalent under homeomorphisms. 
One GOF-knot in $(S^2 \times S^1) \# (S^2 \times S^1)$ is drawn on Figure \ref{fig:ten}. This results from the plumbing of fibered annuli in two $S^2 \times S^1$s. 
Therefore every GOF-knot in $(S^2 \times S^1) \# (S^2 \times S^1)$ is of this position.
\\\\
\\\vspace{0.5in}

\subsection{$L(p,q) \# (S^2 \times S^1) \ (|p| \geq 2)$}
\begin{figure}[htbp]
 \begin{center}
  \includegraphics[width=80mm]{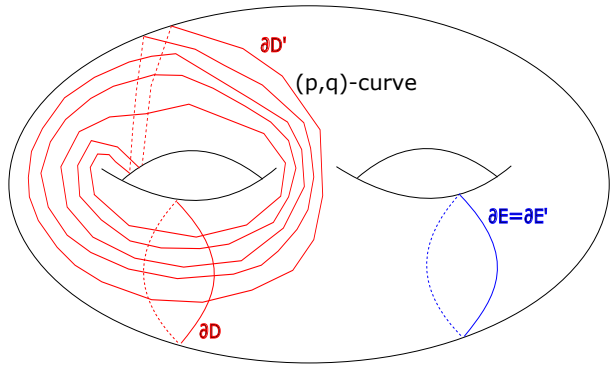}
 \end{center}
 \caption{a standard Heegaard diagram}
 \label{fig:eleven}
\end{figure}

\begin{figure}[htbp]
 \begin{minipage}{0.5\hsize}
  \begin{center}
   \includegraphics[width=50mm]{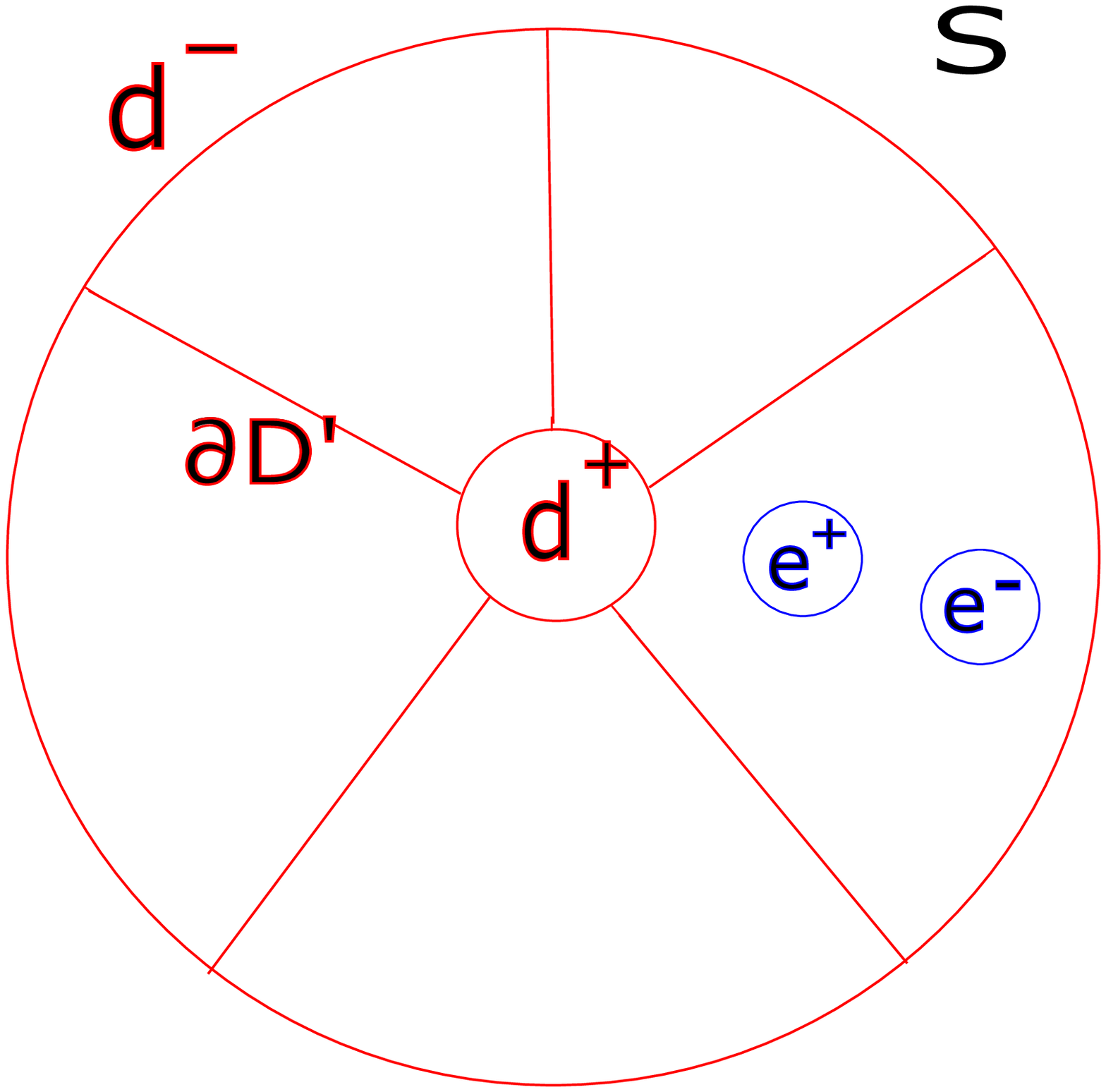}
  \end{center}
  \caption{cut $\Sigma$ along $\partial D$ and $\partial E$}
  \label{fig:twelve}
 \end{minipage}
 \begin{minipage}{0.5\hsize}
  \begin{center}
   \includegraphics[width=30mm]{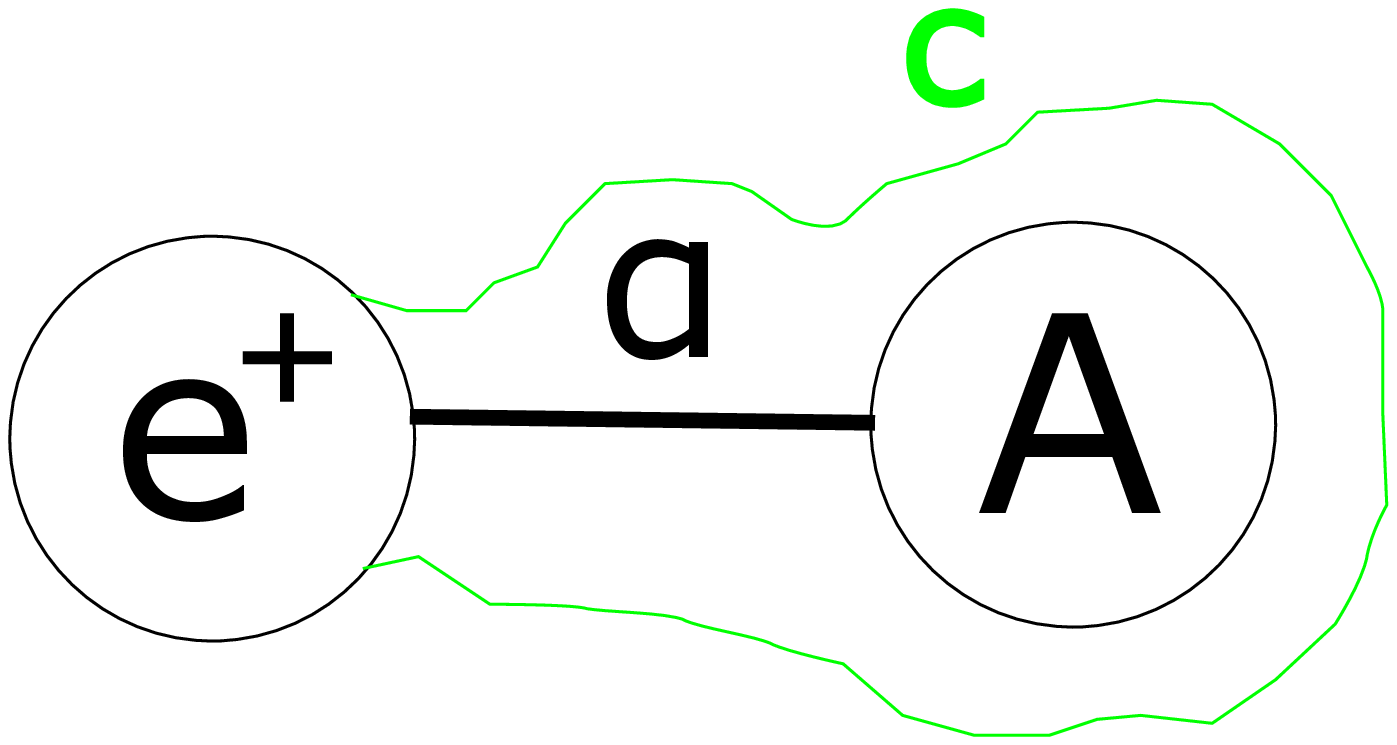}
  \end{center}
  \caption{$yy^{-1}$}
  \label{fig:thirteen}
 \end{minipage}
\end{figure}

\begin{figure}[htbp]
 \begin{center}
  \includegraphics[width=100mm]{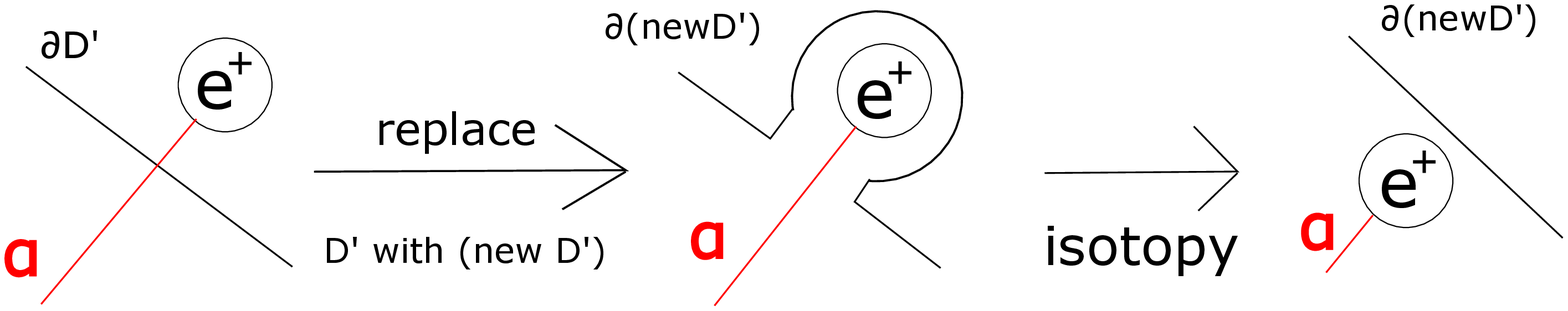}
 \end{center}
 \caption{change $D'$ so that $\alpha$ is disjoint from $\partial D'$}
 \label{fig:fourteen}
\end{figure}

\begin{figure}[htbp]
 \begin{center}
  \includegraphics[width=100mm]{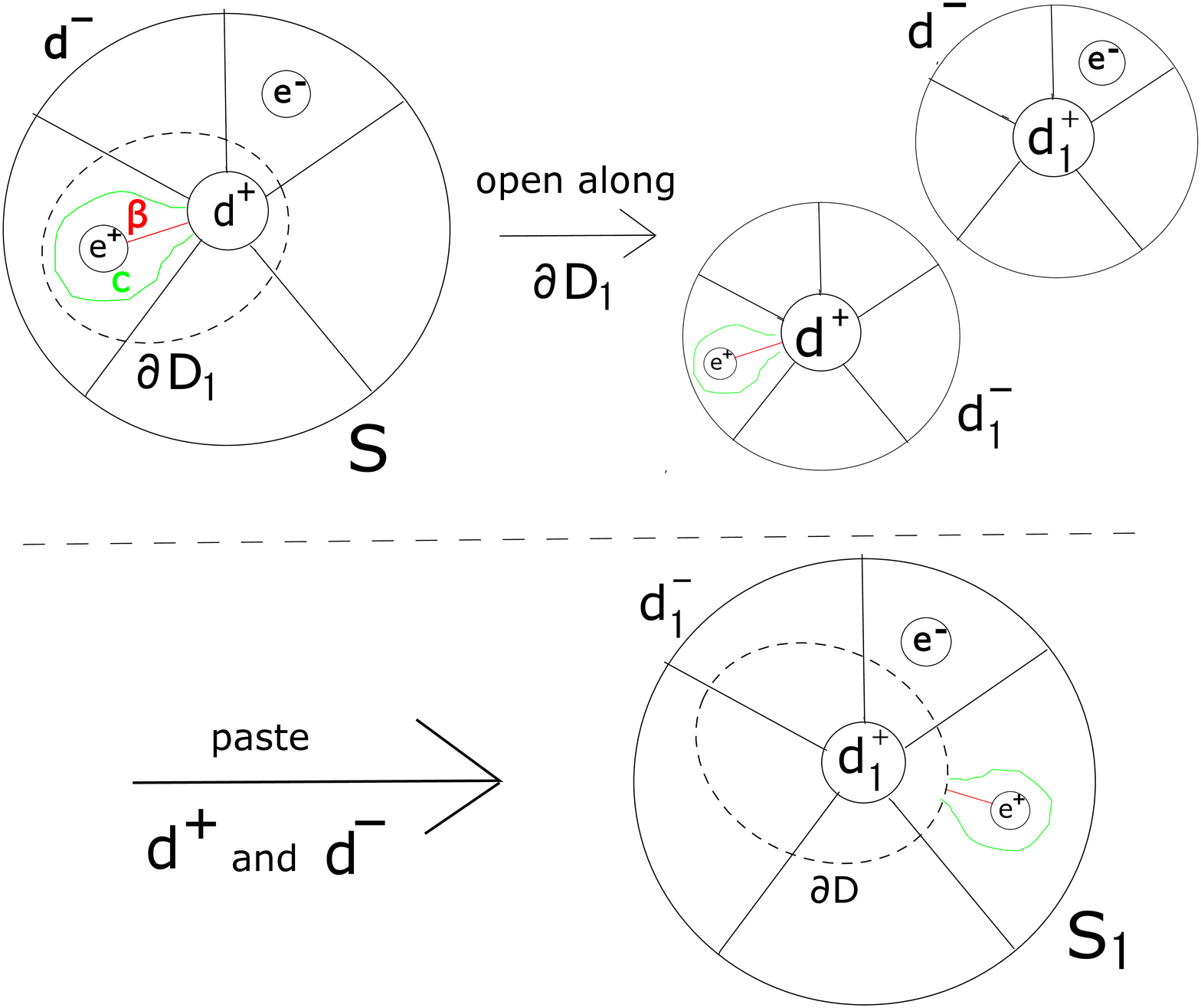}
 \end{center}
 \caption{change $D$ to new disk $D_1$}
 \label{fig:fifteen}
\end{figure}

\begin{figure}[htbp]
 \begin{center}
  \includegraphics[width=60mm]{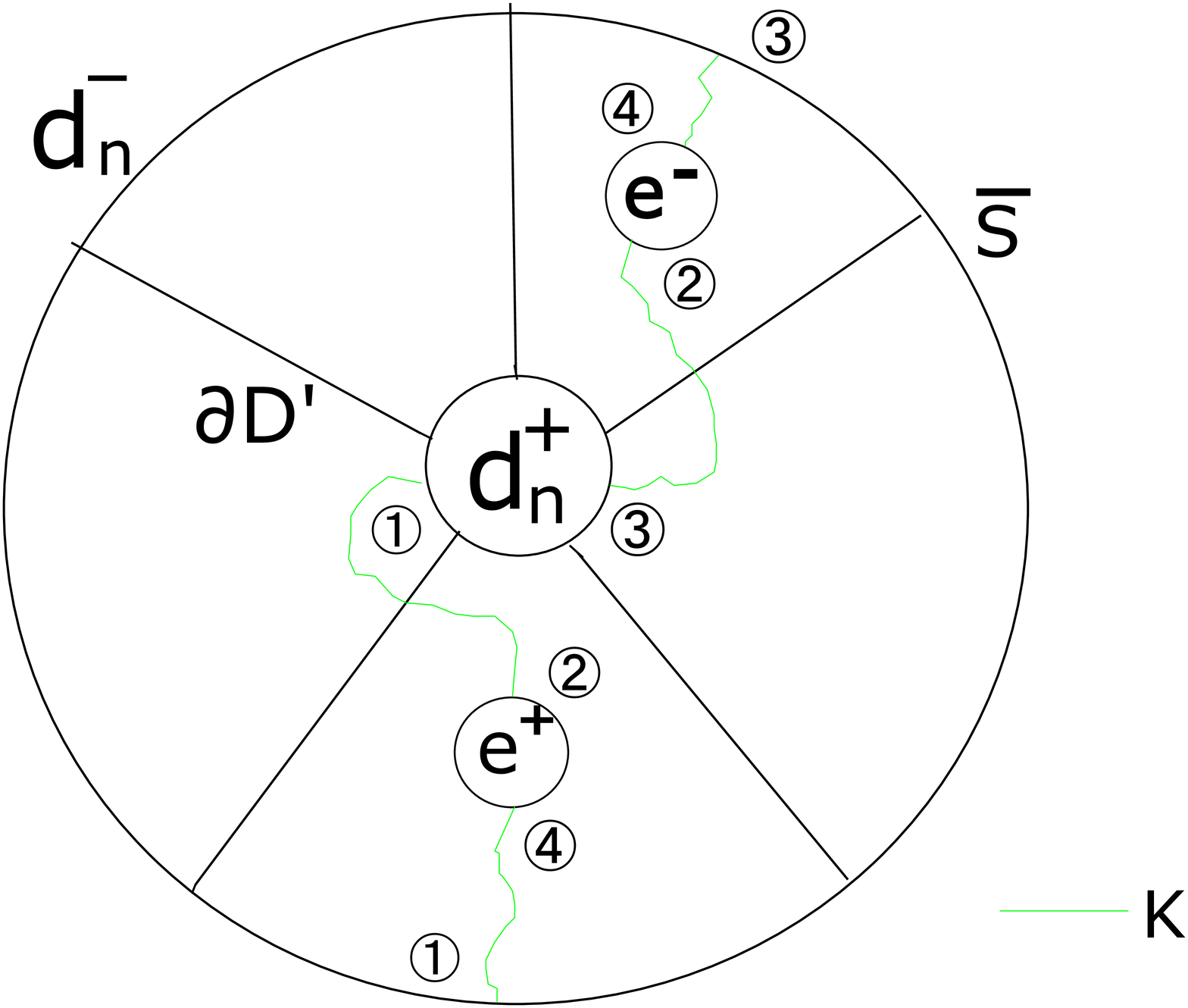}
 \end{center}
 \caption{GOF-knot in $L(p,1)\# (S^2\times S^1)$}
 \label{fig:sixteen}
\end{figure}

\begin{figure}[htbp]
 \begin{center}
  \includegraphics[width=100mm]{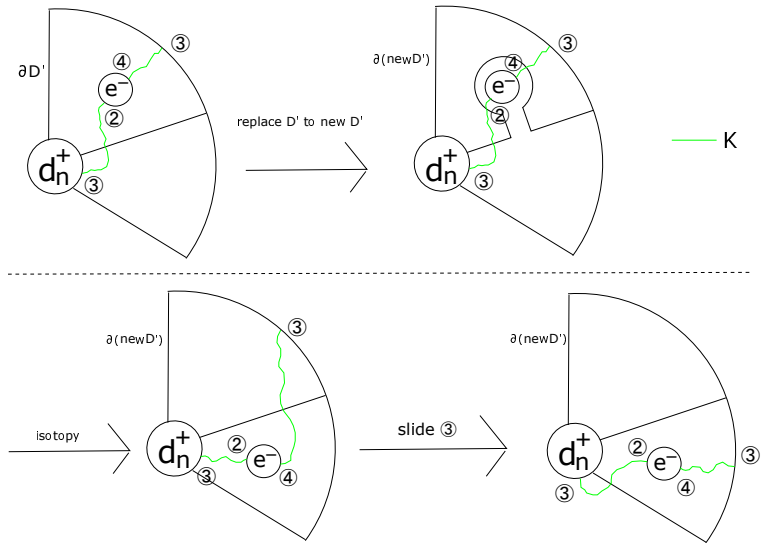}
 \end{center}
 \caption{change cell which contains $e^{-}$}
 \label{fig:seventeen}
\end{figure}

\begin{figure}[htbp]
 \begin{center}
  \includegraphics[width=80mm]{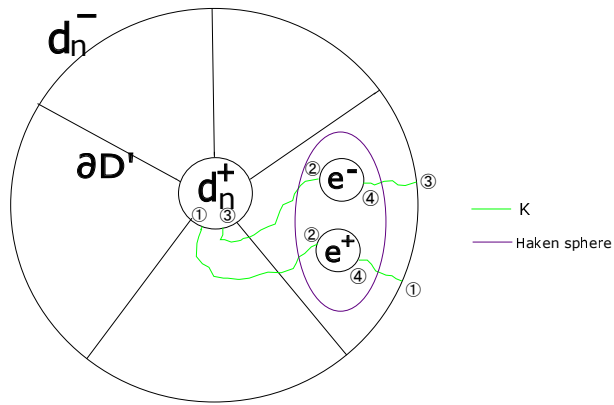}
 \end{center}
 \caption{GOF-knot and Haken sphere in $L(p,1)\# (S^2\times S^1)$}
 \label{fig:eighteen}
\end{figure}

At first, we determine the condition for $L(p,q) \# (S^2 \times S^1)$ to have GOF-knots. Next, we find the positions of GOF-knots if there are. 
In this case, the Heegaard surface is unique as in Section \ref{sec2}. We consider a standard Heegaard diagram in Figure \ref{fig:eleven}. 
We denote by $V \cup _{\Sigma} W$ the corresponding genus two Heegaard splitting. 
We give $\partial D$ and $\partial E$ (resp. $\partial D'$ and $\partial E'$) letters $x$ and $y$ (resp. $x'$ and $y'$). 
We set $S$ to be $\Sigma \setminus (\partial D \cup \partial E)$. It is a sphere with four boundary components. 
We denote by $d^{+}$, $d^{-}$, $e^{+}$ and $e^{-}$ the boundary components of $S$ coming from $\partial D$ and $\partial E$. 
$\partial D'$ cuts $S$ into $p$ cells. See Figure \ref{fig:twelve}. 
In this figure, the subarc of $d^{+}$ which is on the $i$-th cell is identified the subarc of $d^{-}$ which is on $i+q \pmod{p}$. \\
\\
\\
\ We assume there is a GOF-knot $K$. If there is a subarc $c$ of $K$ corresponding to the word of the form $yy^{-1}$, 
there is an arc $\alpha$ in $S$ connecting $e^{+}$ and another boundary component, 
denoted by $A$, such that the boundary of a regular neighborhood of $\alpha \cup A$ in $S$ is $c$ (See Figure \ref{fig:thirteen}). 
If $A$ is $e^{-}$, every subarc of $K$ in $S$ one of whose endpoints is on $e^{-}$ has the other endpoint on $e^{+}$. 
Let $n$ be the number of these subarcs. 
Then the number of the points in $|K \cap e^{+}|$ (in $S$) is at least $(n+2)$. It cannot occur. Therefore $A$ cannot be $e^{-}$. 
So $A$ is $d^{+}$ or $d^{-}$. We can deform the Heegaard diagram so that $\alpha$ is disjoint from $\partial D'$ in the following way: 
If $\alpha$ intersects $\partial D'$, we can get another non-separating disk in $W$ which is disjoint from and not parallel to $D'$ and $E'$ 
by using the subarc of $\alpha$ which connect $e^{+}$ and an intersection point with $\partial D'$ and whose interior is disjoint from $\partial D'$. 
We replace $D'$ with this new disk. Denote it by $D'$ and give the letter $x'$ again. In $S$, the cell which contains $e^{+}$ will change. See Figure \ref{fig:fourteen}. 
Using this operation in finitely many times, we can assume $\alpha$ is disjoint from $\partial D'$. 
In this situation, the subarc $c$ represents a word in $x'$ and $y'$ of the form $y'x'^{\pm p}y'^{-1}$. 
Then by Lemma \ref{lem1} the word in $x'$ and $y'$ represented by $K$ is reduced. Since $|p|\geq 2$, it is not a commutator. 
By Lemma \ref{lem3}, it is a contradiction. 
Hence there are no subarcs of $K$ representing $yy^{-1}$. (Similarly, no $y^{-1}y$.)\\
\\
\ If there is a subarc of $K$ corresponding to a word of the form $xx^{-1}$, there is an arc $\beta$ in $S$ connecting $d^{+}$ and another boundary component, 
denoted by $B$ such that the boundary of a regular neighborhood of $\beta \cup B$ in $S$ is the subarc. 
For the same reason as the above, $B$ is $e^{+}$ or $e^{-}$. 
If $\beta$ intersects $\partial D'$, we replace $D'$ to the new disk which is obtained by disk surgery as in Figure \ref{fig:fourteen} and 
denote this new disk by $D'$ and give the letter $x'$ again.
Using this operation in finitely many times, we can assume $\beta$ is disjoint from $\partial D'$. 
Then as in Figure \ref{fig:fifteen}, we can get a non-separating disk in $V$ by disk surgery. 
Denote this new disk by $D_1$ and give $D_1$ and $E$ the letters $x_1$ and $y_1$. 
We set $S_1$ to be $\Sigma \setminus (\partial D_1 \cup \partial E)$. 
We denote by $d^{+}_1$, $d^{-}_1$, $e^{+}$ and $e^{-}$ the boundary components of $S_1$ coming from $\partial D_1$ and $\partial E$. 
In $S_1$, there cannot be subarcs of $K$ of the form $y_1y_{1}^{-1}$ by the same discussion as the above.\\
\ Inductively, if there is a subarc of $K$ of the form $x_kx_{k}^{-1}$, there is an arc $\beta _k$ in $S_k$ connecting $d_{k}^{+}$ and another boundary component, 
denoted by $B_k$ such that the boundary of a regular neighborhood of $\beta _k \cup B_k$ in $S_k$ is the subarc. 
For the same reason as the above, $B_k$ is $e^{+}$ or $e^{-}$. 
If $\beta _k$ intersects $\partial D'$, 
we replace $D'$ to the new disk which is obtained by disk surgery and denote this new disk by $D'$ and give the letter $x'$ again.
Using this operation in finitely many times, we can assume $\beta_k$ is disjoint from $\partial D'$. Then as in Figure \ref{fig:fifteen}, 
we can get a non-separating disk in $V$ by disk surgery. Denote this new disk by $D_{k+1}$ and give $D_{k+1}$ and $E$ the letters $x_{k+1}$ and $y_{k+1}$. 
We set $S_{k+1}$ to be $\Sigma \setminus (\partial D_{k+1} \cup \partial E)$. We denote by $d^{+}_{k+1}$, $d^{-}_{k+1}$, $e^{+}$ and $e^{-}$ 
the boundary components of $S_{k+1}$ coming from $\partial D_{k+1}$ and $\partial E$. 
In $S_{k+1}$, there cannot be subarcs of $K$ of the form $y_{k+1}y_{k+1}^{-1}$ by the same discussion as the above. 
Then for some non-negative integer $n$, the word in {$x_n$, $y_n$} represented by $K$ is reduced. We regard $x_0,y_0$ and $S_0$ as $x,y$ and $S$. 
In $S_n$, $K$ is a collection of four arcs, 
connecting $e^{+}$ and $d^{+}_{n}$, connecting $d^{-}_{n}$ and $e^{+}$, connecting $e^{-}$ and $d^{-}_{n}$ and connecting ${d_{n}}^{+}$ and $e^{-}$ 
(the terminal point of an arc is the starting point of the next arc.). 
In this situation, there cannot be a subarc of $K$ of the form $y'y'^{-1}$ and by changing $D'$ to a new disk 
(and denoting this new disk $D'$ and giving the letter $x'$ again), we can assume there are no subarcs of $K$ of the form $x'x'^{-1}$ (See Figure \ref{fig:fourteen}). 
Let $\bar{S}$ denote a sphere with four boundary components obtained by $S_n$ and changing $\partial D'$ as in above (See Figure \ref{fig:sixteen}). 
In $\bar{S}$, $K$ is represented by four arcs and is simultaneously a reduced form in $\{x',y'\}$. 
Hence it is necessary that $q \equiv \pm 1$  mod $p$ since the word represented by $K$ in $\{x',y'\}$ contains ${x'}^{[q]}$ or ${x'}^{[p-q]}$. 
($[n]$ is the residue class of $n$ mod $p$.) $K$ is of the position in Figure \ref{fig:sixteen} for example.\\
\\
\ In $\bar{S}$, by repeating disk surgeries as in Figure \ref{fig:seventeen}, we can assume that $e^{+}$ and $e^{-}$ are in the same cell. 
By looking the Haken sphere in Figure \ref{fig:eighteen}, we see that $K$ with its fiber $T$ is the result of the plumbing of two fibered annuli 
in $L(p,q)$ and $S^2 \times S^1$ respectively. 
By changing the orientation if necessary, we assume that $L(p,q)$ is $L(p,1)$ and the fibered annulus in it is the $p$-Hopf band. 
Then, this GOF-knot is unique under self-homeomorphisms.\\
\\
\\
\ Therefore, we conclude that the necessary and sufficient condition for $L(p,q) \# (S^2 \times S^1)$ to have GOF-knots is $q \equiv \pm 1$  mod $p$, 
and if there is a GOF-knot, it is unique and obtained by the plumbing. 
\\\vspace{0.5in}

\subsection{$L(p_1,q_1) \# L(p_2,q_2) \ (|p_1|,|p_2| \geq 2)$}
\begin{figure}[htbp]
 \begin{center}
  \includegraphics[width=75mm]{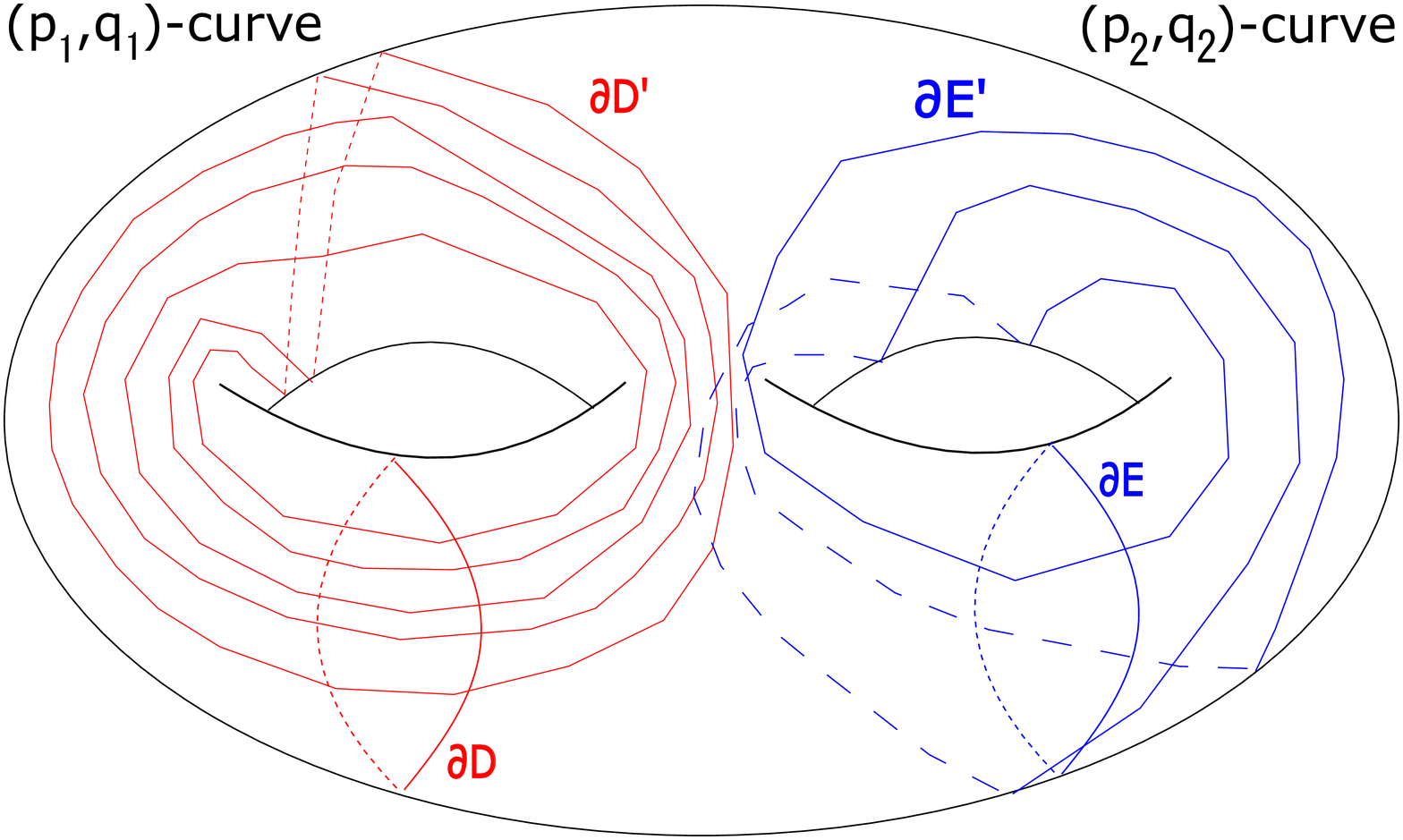}
 \end{center}
 \caption{a standard Heegaard diagram}
 \label{fig:nineteen}
\end{figure}

\begin{figure}[htbp]
 \begin{center}
  \includegraphics[width=45mm]{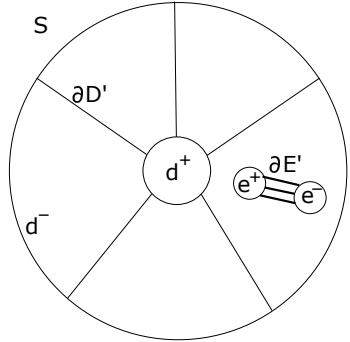}
 \end{center}
 \caption{cut $\Sigma$ along $\partial D$ and $\partial E$}
 \label{fig:twenty}
\end{figure}

\begin{figure}[htbp]
 \begin{center}
  \includegraphics[width=80mm]{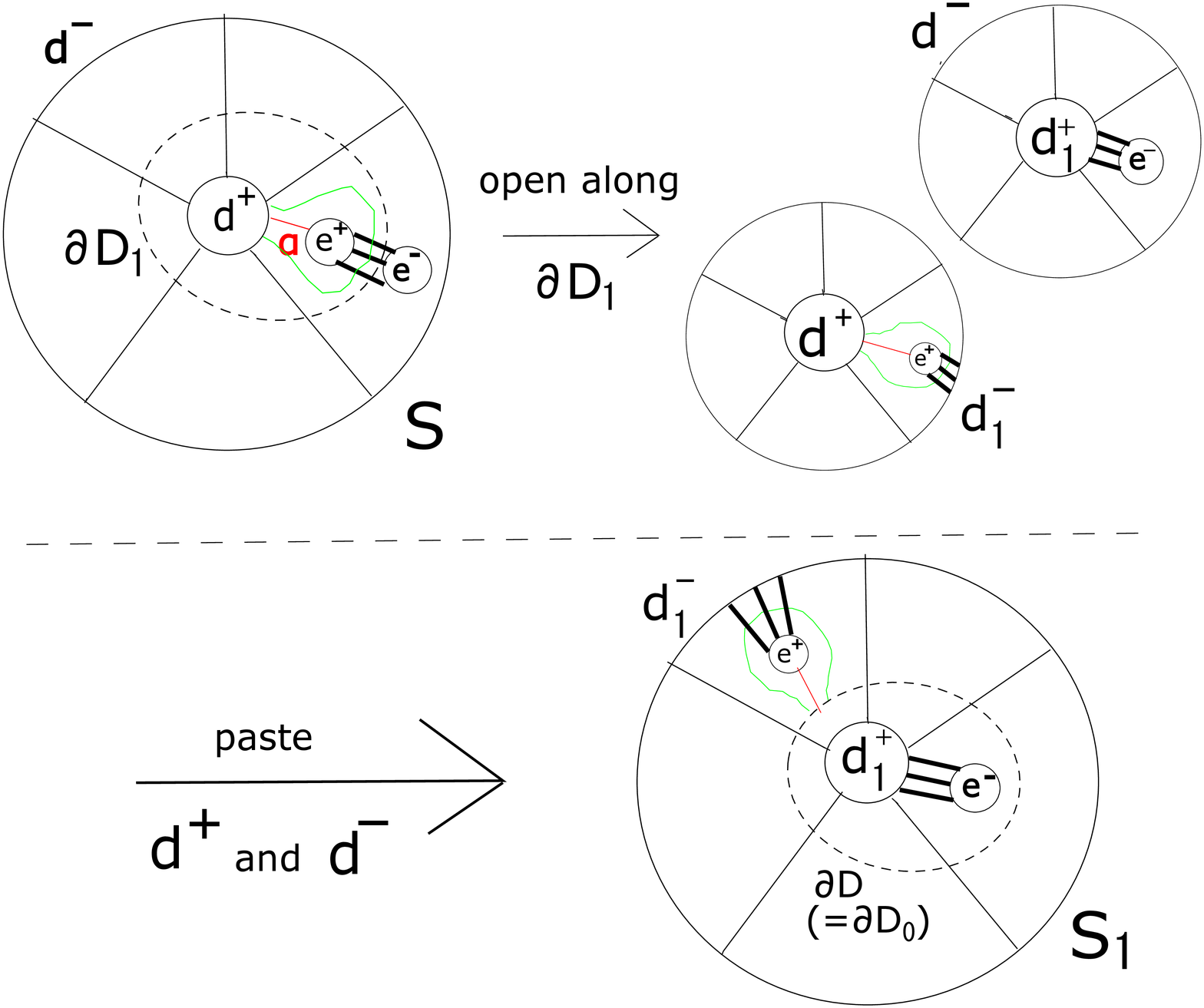}
 \end{center}
 \caption{make $S_1$ from $S$(=$S_0$)}
 \label{fig:twentyone}
\end{figure}

\begin{figure}[htbp]
 \begin{center}
  \includegraphics[width=100mm]{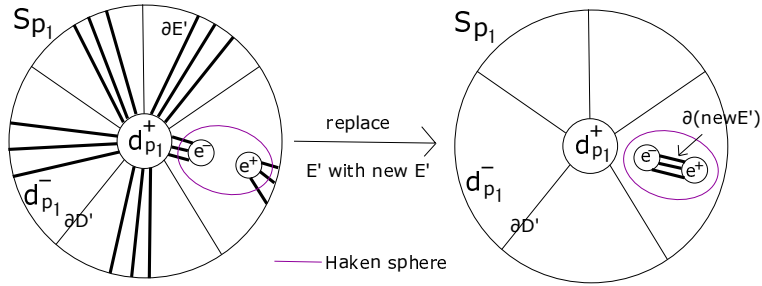}
 \end{center}
 \caption{another standard Heegaard diagram}
 \label{fig:twentytwo}
\end{figure}

\begin{figure}[htbp]
 \begin{center}
  \includegraphics[width=30mm]{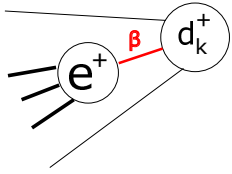}
 \end{center}
 \caption{$\beta$ connects $e^{+}$ and $d_k^{\epsilon}$. $\beta$ is not connected directly by $\partial E'$}
 \label{fig:twentythree}
\end{figure}

\begin{figure}[htbp]
 \begin{center}
  \includegraphics[width=100mm]{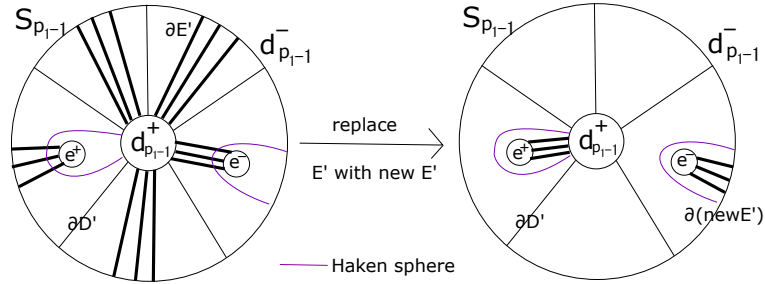}
 \end{center}
 \caption{another Heegaard diagram like $S_1$}
 \label{fig:twentyfour}
\end{figure}

\begin{figure}[htbp]
 \begin{center}
  \includegraphics[width=100mm]{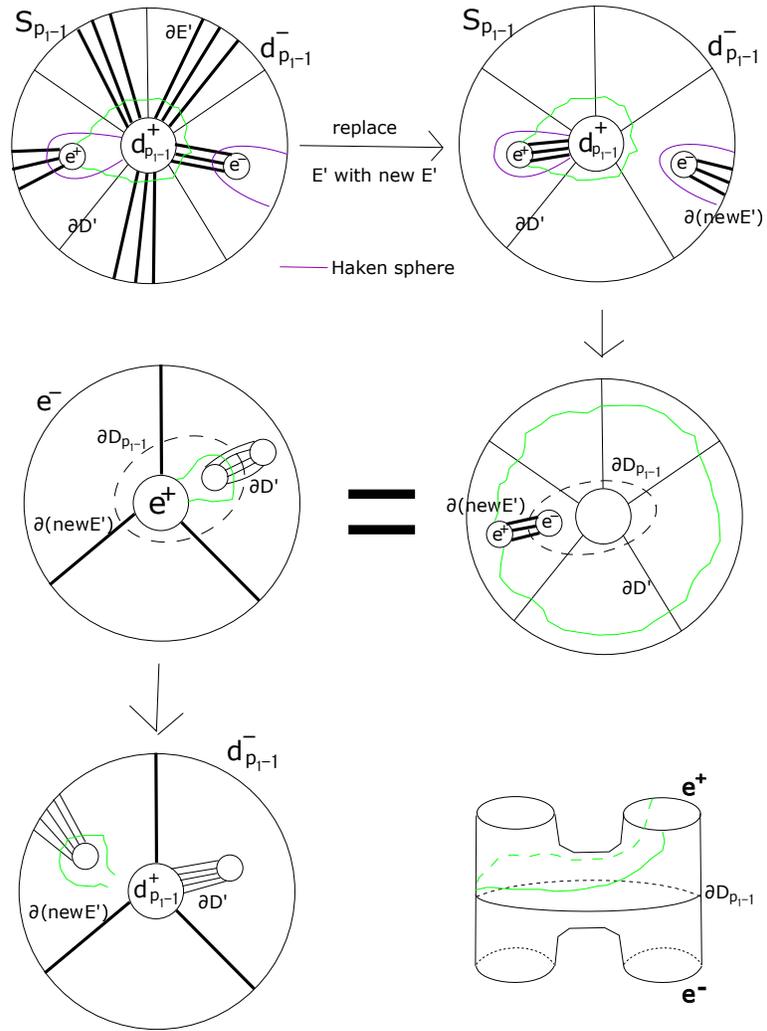}
 \end{center}
 \caption{the case having $y_{p_1-1}y_{p_1-1}^{-1}$}
 \label{fig:twentyfive}
\end{figure}

\begin{figure}[htbp]
 \begin{center}
  \includegraphics[width=100mm]{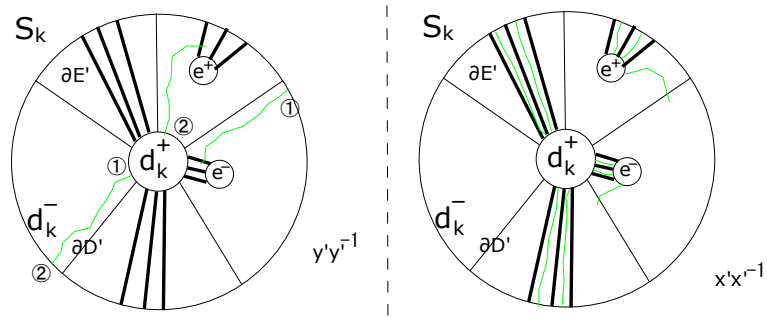}
 \end{center}
 \caption{subarc of type $y'y'^{-1}$ and $x'x'^{-1}$ in $S_k$}
 \label{fig:twentysix}
\end{figure}

\begin{figure}[htbp]
 \begin{center}
  \includegraphics[width=80mm]{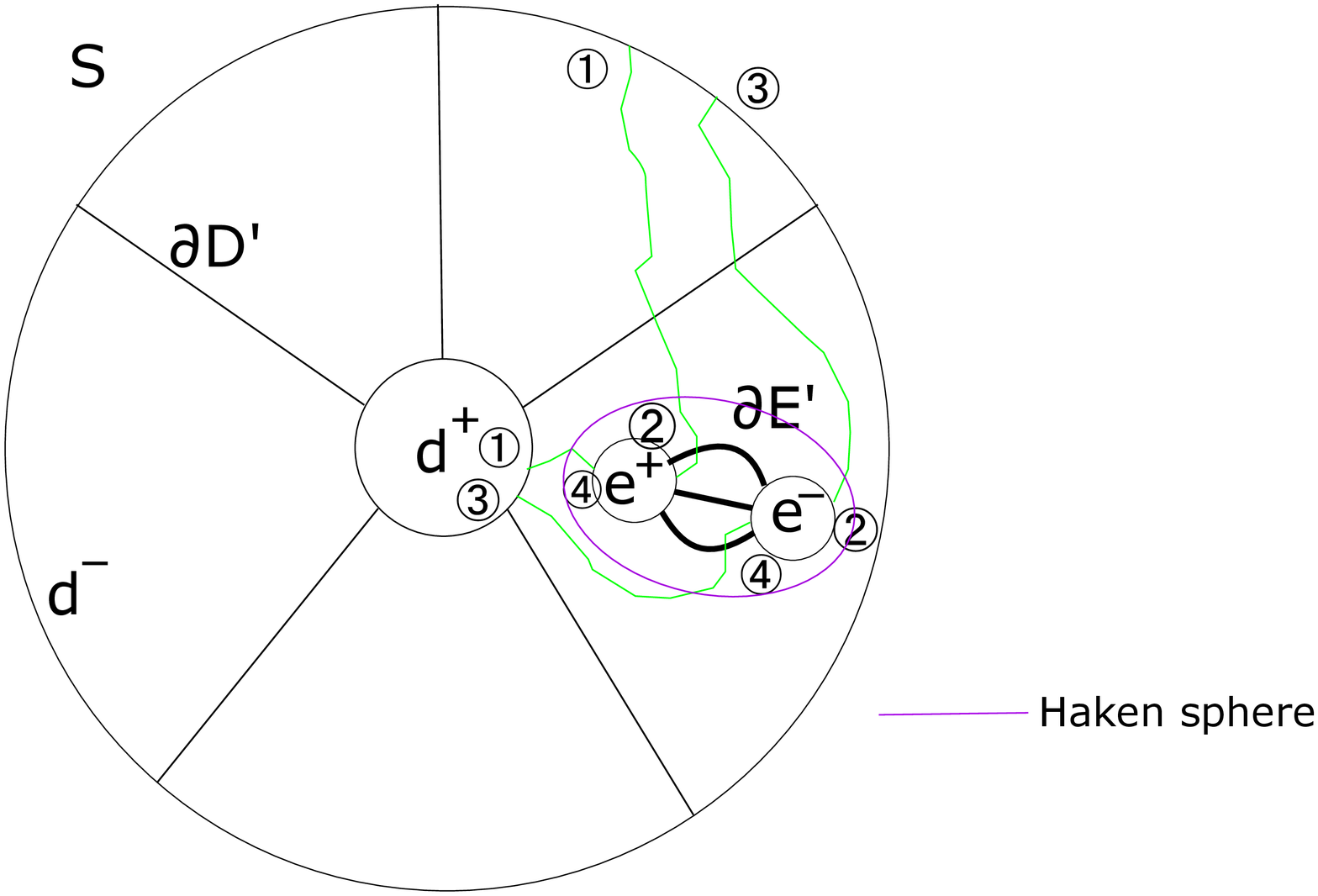}
 \end{center}
 \caption{$K$ in a standard Heegaard splitting}
 \label{fig:twentyseven}
\end{figure}

At first, we determine the condition for $L(p_1,q_1) \# L(p_2,q_2)$ to have GOF-knots. Next, we find the positions of GOF-knots if there are. 
In this case, there are at most two Heegaard surfaces of genus two as in Section \ref{sec2}. However the corresponding Heegaard diagrams of them are similar. 
Thus we assume that a GOF-knot $K$ is on a standard Heegaard surface in Figure \ref{fig:nineteen}. We denote by $V \cup _{\Sigma} W$ this Heegaard splitting. 
We give $\partial D$ and $\partial E$ (resp. $\partial D'$ and $\partial E'$) letters $x$ and $y$ (resp. $x'$ and $y'$). 
We set $S$ to be $\Sigma \setminus (\partial D \cup \partial E)$. It is a sphere with four boundary components. 
We denote by $d^{+}$, $d^{-}$, $e^{+}$ and $e^{-}$ the boundary components of $S$ coming from $\partial D$ and $\partial E$ (See Figure \ref{fig:twenty}). 
\\
\\
\ Since $xx^{-1}$ and $yy^{-1}$ cannot coexist by Lemma \ref{rem}, we can assume that there are no subarcs of type $yy^{-1}$. 
If there is a subarc $c$ of $K$ representing $xx^{-1}$ (by the symmetry, same for $yy^{-1}$), 
there is an arc $\alpha$ on $S$ connecting $d^{+}$ and $e^{\epsilon}$ ($\epsilon \in \{ +,- \}$) 
such that the boundary of a regular neighborhood of $\alpha \cup e^{\epsilon}$ in $S$ is $c$. 
We say such an arc $\alpha$ is the corresponding arc of $c$. 
We assume $\epsilon$ is $+$. Isotoping $\alpha$ so that it intersects $\partial D'$ and $\partial E'$ minimally and by Lemma \ref{lem1}, 
we can assume that $\alpha$ is disjoint from $\partial D'$ and $\partial E'$. 
This is because if $\alpha$ intersects them and if the intersection point nearest to $e^{+}$ is on $\partial E'$, 
this intersection point can be omitted by isotopy and if the intersection point nearest to $e^{+}$ is on $\partial D'$, 
the subarc of $c$ (so of $K$) represents a word of type $x'{y'}^{p_2}{x'}^{-1}$ and then by Lemma \ref{lem1}, 
$K$ is a reduced form in $\{x',y'\}$ and has $y'^{p_2}$ so $K$ cannot be the commutator of $x'$ and $y'$ after reduction. 
We set $D_0$ to be $D$ and set $D_1$ to be the non-separating disk in $V$ obtained by disk surgery of $D_0$ using $c$ (See Figure \ref{fig:twentyone}). 
In Figure \ref{fig:twentyone}, $S_1$ is the sphere with four boundary components obtained by cutting $\Sigma$ along $\partial D_1$ and $\partial E$. 
We denote by $d_{1}^{+}$, $d_{1}^{-}$, $e^{+}$ and $e^{-}$ the boundary components of $S_1$ coming from $\partial D_1$ and $\partial E$. 
We give $\partial D_1$ and $\partial E$ the letters $x_1$ and $y_1$. \\
\\
\ If $K$ has a subarc representing a word of type $x_1x_1^{-1}$, there is a corresponding arc on $S_1$ which connects $d_1^{+}$ and $e^{+}$. 
Note that this corresponding arc does not connect $d_1^{+}$ and $e^{-}$. 
By an argument similar to the above, we can assume this corresponding arc is disjoint from $\partial D'$ and $\partial E'$. 
Inductively for $k$ ($0 \leq k \leq p_1-1$), if $K$ has a subarc representing a word of type $x_kx_k^{-1}$ there is a corresponding arc 
which connects $d_k^{+}$ and $e^{+}$. 
We can assume this corresponding arc is disjoint from $\partial D'$ and $\partial E'$. 
We set $D_{k+1}$ to be the non-separating disk in $V$ obtained by the disk surgery of $D_k$ using the subarc of $K$. 
Let $S_{k+1}$ be the sphere with four boundary components obtained by cutting $\Sigma$ along $\partial D_{k+1}$ and $\partial E$. 
We denote by $d_{k+1}^{+}$, $d_{k+1}^{-}$, $e^{+}$ and $e^{-}$ the boundary components of $S_{k+1}$ coming from $\partial D_{k+1}$ and $\partial E$. 
We give $\partial D_{k+1}$ and $\partial E$ the letters $x_{k+1}$ and $y_{k+1}$. (We regard $S$, $D$ and $d^{\pm}$ as $S_0$, $D_0$ and $d_{0}^{\pm}$.) 
If this operation is repeated until we get $S_{p_1}$, we can get another standard Heegaard diagram (See Figure \ref{fig:twentytwo}). 
In this standard Heegaard diagram, the intersection number of $K$ and this new pair of disjoint non-separating, non-parallel disks in $V$ is less than 
that of $K$ with old pair of disjoint non-separating, non-parallel disks in $V$. 
Thus we can assume this operation stops at getting $S_k$ ($0 \leq k \leq p_1-1$).\\
\\
\ Note that if there is a subarc of $K$ of type $y_ky_k^{-1}$ ($1\leq k\leq p_1-1$) in $S_k$, its corresponding arc $\beta$ is like in Figure \ref{fig:twentythree}. 
Otherwise there is a subarc of $K$ of type $y_{k-1}y_{k-1}^{-1}$ in $S_{k-1}$. 
It contradicts Lemma \ref{rem}, $x_{k-1}x_{k-1}^{-1}$ and $y_{k-1}y_{k-1}^{-1}$ cannot coexist.\\
\\
\ Moreover, if this operation is repeated until we get $S_{p_1-1}$ and there are no $x_{p_1-1}x_{p_1-1}^{-1}$, 
we can get another Heegaard diagram which is obtained by the standard diagram of the above. It is similar to $S_1$ (See Figure \ref{fig:twentyfour}). 
If there is a subarc of $K$ of type $y_{p_1-1}y_{p_1-1}^{-1}$, we can get another Heegaard splitting like in Figure \ref{fig:twentyfive}. In this Heegaard diagram, 
the intersection number of the pair of disjoint non-separating, non-parallel disks in $V$ with $K$ is not more than that of $D_{p_1-1} \cup E$ and $K$.
Thus we can assume this operation stops at getting $S_k$ ($0 \leq k \leq p_1-2$).\\
\\
\ In this setting, if there is subarc of $K$ representing $y_ky_k^{-1}$ (or $y_k^{-1}y_k$) ($k\neq 0$ from our assumption.), 
there are no subarcs of $K$ representing the words of type $x'x'^{-1}$ nor $y'y'^{-1}$ (See Figure \ref{fig:twentysix}). 
If there are no subarcs of $K$ representing $y_ky_k^{-1}$ (or $y_k^{-1}y$), 
there are no subarcs of $K$ representing the words of type $x'x'^{-1}$ nor $y'y'^{-1}$ neither (this is because $K$ is a reduced form in $S_k$.) 
(See Figure \ref{fig:twentysix}). 
Thus $K$ is a reduced form in $D'$ and $E'$. Moreover in this situation $K$ must be a reduced form in $D_0$ and $E_0$ too. 
This is because a subarc in $S_0$ representing $x_0x_0^{-1}$ in $\{x_0,y_0\}$ represents $y'^{\pm p_2}$ in $\{x',y'\}$. 
\\
\ As a result, $K$ must be a reduced form in a standard Heegaard diagram (in Figure \ref{fig:twentyseven}) in $\{x,y\}$ and $\{x',y'\}$ simultaneously. 
Hence it is necessary that $q_i \equiv \pm 1$  mod $p_i$ ($i = 1,2$) as in 4.2 , 
and $K$ is of the position in Figure \ref{fig:twentyseven} for example.\\
\\
\ By looking the Haken sphere in Figure \ref{fig:twentyseven}, we see that $K$ with its fiber $T$ is the result of the plumbing of two fibered annuli 
in $L(p_1,q_1)$ and $L(p_2,q_2)$ respectively. This implies $q_i \equiv \pm 1$  mod $p_i$ ($i=1,2$), 
and then the condition to have the unique genus two Heegaard splitting is satisfied.\\
\ By changing the orientation if necessary, we assume $L(p_1,q_1)$ is $L(p_1,1)$ and the fibered annulus in it is $p_1$-Hopf band. 
Note that as noted in \cite{10}, $L(r_1,s_1)\# L(r_2,s_2)$ is homeomorphic to $L(r_1,s_1)\# L(r_2,-s_2)$ 
if and only if ${s_1}^{2}\equiv -1$ mod $r_1$ or ${s_2}^{2}\equiv -1$ mod $r_2$. 
Hence if neither $p_1$ nor $p_2$ is $2$, the GOF-knot is unique. \\
\ If either $p_1$ or $p_2$ is $2$ (we assume $p_2$ is $2$), there can be two GOF-knots, 
one is obtained by the plumbing of the $p_1$-Hopf band in $L(p_1,1)$ and the $2$-Hopf band in $L(2,1)$ and 
the other is obtained by the plumbing of the $p_1$-Hopf band in $L(p_1,1)$ and the $-2$-Hopf band in $L(2,1)$. 
The monodromy of the former is represented in $GL_{2}(\mathbb{Z})$ as $ \left(\begin{array}{ccc} 1 & p_1 \\2 & 1+2p_1 \end{array} \right)$ and 
that of the latter is represented in $GL_{2}(\mathbb{Z})$ as $ \left(\begin{array}{ccc} 1 & p_1 \\-2 & 1-2p_1 \end{array} \right)$. 
Since they are not conjugate in $GL_{2}(\mathbb{Z})$, these GOF-knots are not equivalent .\\\\

\ Therefore, we conclude that the necessary and sufficient condition for $L(p_1,q_1) \# L(p_2,q_2)$ to have GOF-knots is 
$q_i \equiv \pm 1$  mod $p_i$ ($i=1,2$), and if there is a GOF-knot, it is obtained by the plumbing. 
Moreover, if neither $p_1$ nor $p_2$ is $2$, the GOF-knot is unique and otherwise there are just two GOF-knots. 
\\\vspace{0.5in}

\subsection{$S^2 \times S^1$}
\begin{figure}[htbp]
 \begin{center}
  \includegraphics[width=60mm]{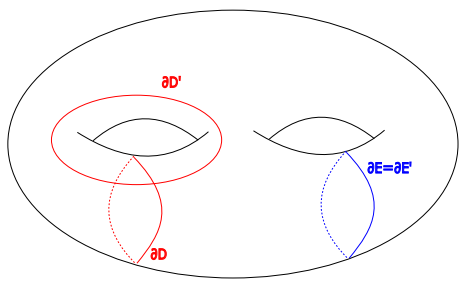}
 \end{center}
 \caption{a standard Heegaard diagram}
 \label{fig:twentyeight}
\end{figure}

\begin{figure}[htbp]
 \begin{center}
  \includegraphics[width=40mm]{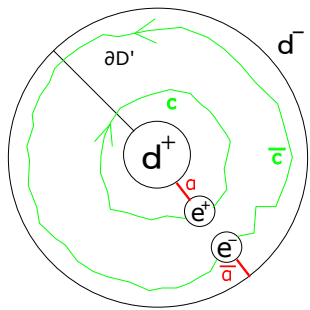}
 \end{center}
 \caption{$c$ and $\bar{c}$}
 \label{fig:twentynine}
\end{figure}

\begin{figure}[htbp]
 \begin{center}
  \includegraphics[width=80mm]{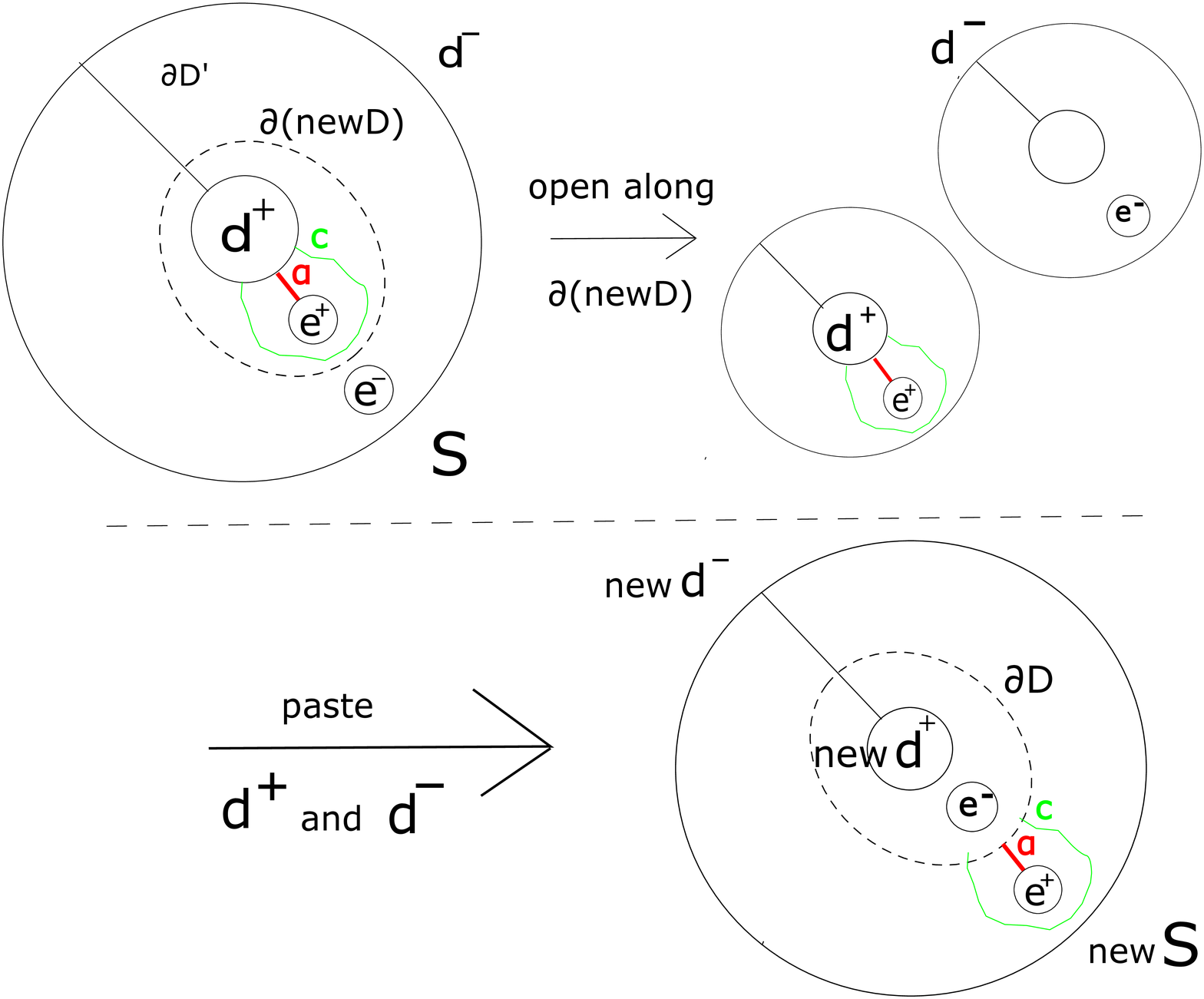}
 \end{center}
 \caption{replace $D$ with new $D$}
 \label{fig:thirty}
\end{figure}

\begin{figure}[htbp]
 \begin{center}
  \includegraphics[width=80mm]{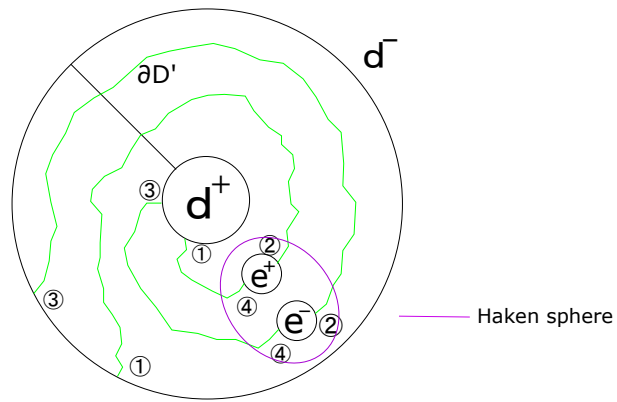}
 \end{center}
 \caption{$K$ on a standard Heegaard splitting}
 \label{fig:thirtyone}
\end{figure}

In this case the genus two Heegaard surface is unique as in Section \ref{sec2}. We regard $S^2 \times S^1$ as $S^3 \# (S^2 \times S^1)$ and 
we consider a standard Heegaard diagram in Figure \ref{fig:twentyeight}. We denote by $V \cup _{\Sigma} W$ the corresponding genus two Heegaard splitting. 
We give $\partial D$ and $\partial E$ (resp. $\partial D'$ and $\partial E'$) letters $x$ and $y$ (resp. $x'$ and $y'$). 
We set $S$ to be $\Sigma \setminus (\partial D \cup \partial E)$. It is a sphere with four boundary components. 
We denote by $d^{+}$, $d^{-}$, $e^{+}$ and $e^{-}$ the boundary components of $S$ coming from $\partial D$ and $\partial E$. \\
\\
\ The argument is almost similar to that of the case $L(p,q) \# (S^2 \times S^1)$. We assume there is a GOF-knot $K$. 
If there is a subarc $c$ of $K$ representing $yy^{-1}$ (or $y^{-1}y$), we can assume the corresponding arc $\alpha$ is disjoint from $\partial D'$. 
In this situation, $c$ represents a word $y'xy'^{-1}$ in $\{x',y'\}$ and by Lemma \ref{lem1} the word in $\{x',y'\}$ represented by $K$ is reduced. 
Moreover, by Lemma \ref{rem} there is a subarc $\bar{c}$ of $K$ representing $y^{-1}y$ and 
we can also assume its representing arc $\bar{\alpha}$ is disjoint from $\partial D'$ as 4.3 . See Figure \ref{fig:twentynine}. 
$\partial D'$ intersects each of $c$ and $\bar{c}$ once and $\partial E'$ intersects each of $c$ and $\bar{c}$ twice. 
Since the word in $\{x',y'\}$ represented by $K$ is reduced, $K$ is $c \cup \bar{c}$. 
In this situation, since $K$ is disjoint from $D$, the word in $\{x,y\}$ represented by $K$ cannot contain $x$. 
Especially $K$ is not a GOF-knot. It is a contradiction. Hence there are no subarcs representing $yy^{-1}$.\\
\\
\ If there is a subarc $c$ of $K$ representing $xx^{-1}$ (or $x^{-1}x$), we assume its corresponding arc $\alpha$ is disjoint from $\partial D'$ as 4.3 . 
In this situation, we replace $D$ with new $D$ as in Figure \ref{fig:thirty} and give new $D$ the letters $x$ and $E$ the letters $y$. 
In this new $S$, there are no subarcs of $K$ representing $yy^{-1}$ by the same discussion above. 
By iterating this operation in finitely many times, $K$ is a reduced form in $\{x,y\}$ in Figure \ref{fig:thirtyone}. 
For the same reason, we can assume $K$ is a reduced form in $\{x',y'\}$ too.\\
\\
\\
\ $K$, which is a reduced form in $\{x,y\}$ and $\{x',y'\}$, is like in Figure \ref{fig:thirtyone}. 
It can be decomposed into two fibered annuli in $S^3$ and $S^2 \times S^1$. 
There can be two GOF-knots, one is obtained by the plumbing of the $+1$-Hopf annulus in $S^3$ and the fibered annulus in $S^2\times S^1$ and 
the other is obtained by the plumbing of the $-1$-Hopf annulus in $S^3$ and the fibered annulus in $S^2\times S^1$. 
The monodromy of the former (it is the positive Dehn twist along the simple closed curve on a fiber corresponding to the core curve of $+1$-Hopf annulus) 
is represented in $GL_{2}(\mathbb{Z})$ as $ \left(\begin{array}{ccc} 1 & 1 \\0 & 1 \end{array} \right)$ and 
that of the latter (it is the negative Dehn twist along the simple closed curve on a fiber corresponding to the core curve of $-1$-Hopf annulus) 
is represented in $GL_{2}(\mathbb{Z})$ as $ \left(\begin{array}{ccc} 1 & -1 \\0 & 1 \end{array} \right)$. 
Since they are conjugate in $GL_{2}(\mathbb{Z})$, these GOF-knots are equivalent.\\\\
\ Therefore, we conclude that $S^2 \times S^1$ has the unique GOF-knot and it (and its fiber) are obtained by the plumbing of two fibered annuli in $S^3$ and $S^2 \times S^1$.
\\\vspace{0.5in}

\subsection{$S^3$}
\begin{figure}[htbp]
 \begin{center}
  \includegraphics[width=60mm]{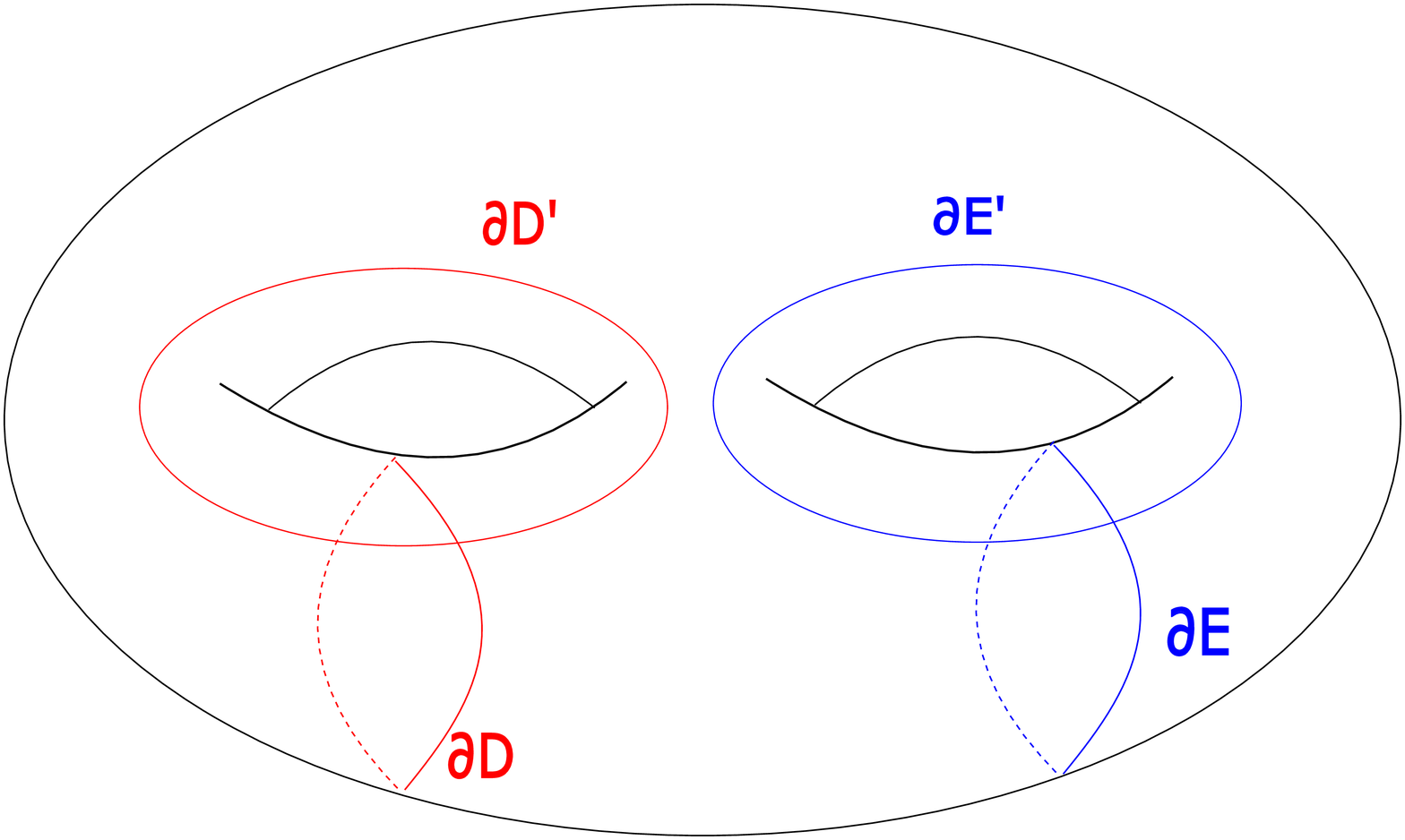}
 \end{center}
 \caption{a standard Heegaard diagram}
 \label{fig:thirtytwo}
\end{figure}

\begin{figure}[htbp]
 \begin{center}
  \includegraphics[width=100mm]{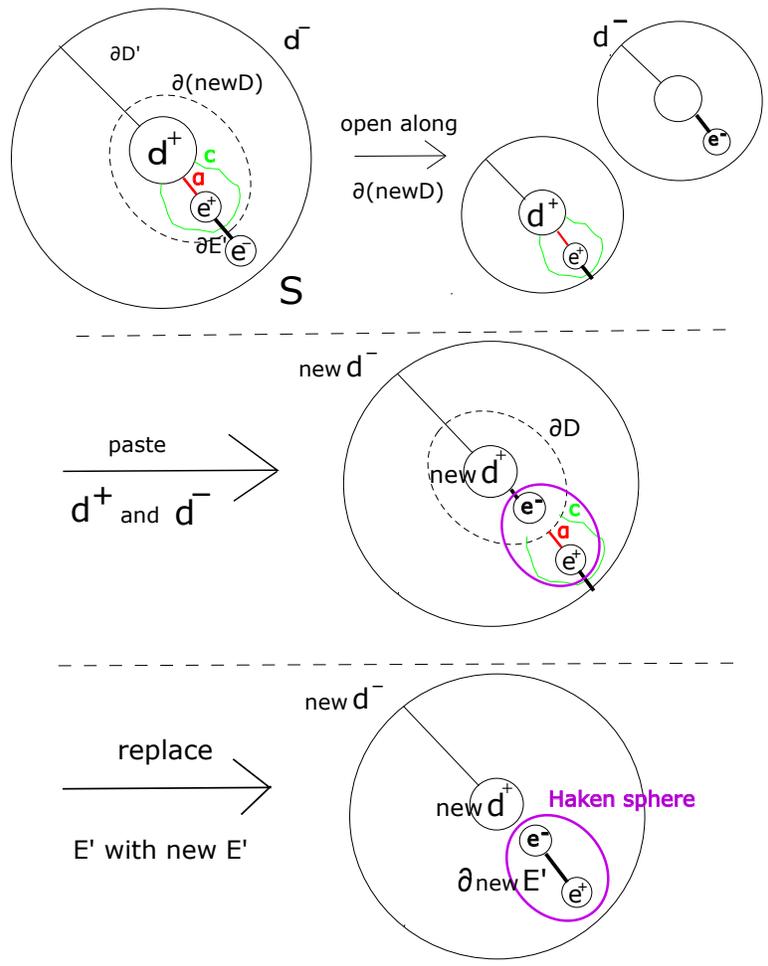}
 \end{center}
 \caption{another standard Heegaard diagram}
 \label{fig:thirtythree}
\end{figure}

\begin{figure}[htbp]
 \begin{center}
  \includegraphics[width=100mm]{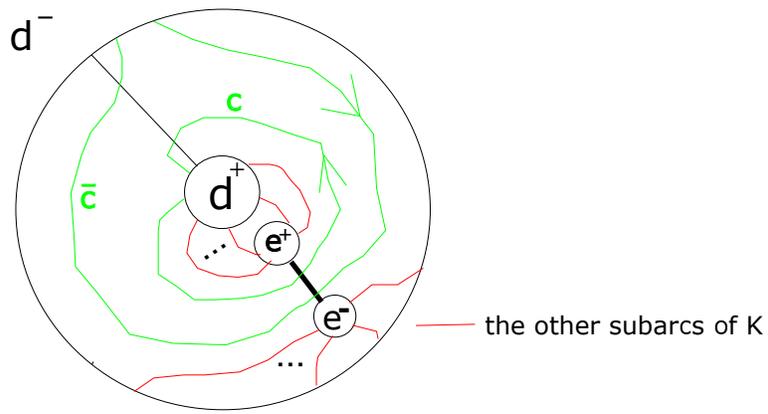}
 \end{center}
 \caption{$c$, $\bar{c}$ and the other subarcs of $K$}
 \label{fig:thirtyfour}
\end{figure}

\begin{figure}[htbp]
 \begin{center}
  \includegraphics[width=120mm]{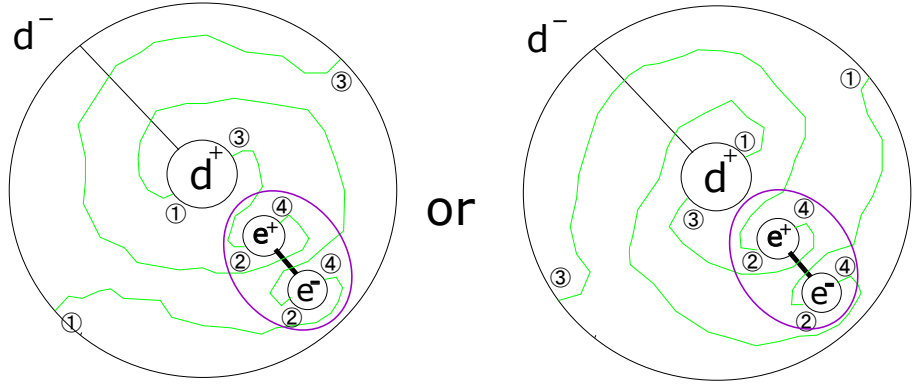}
 \end{center}
 \caption{reduced $K$}
 \label{fig:thirtyfive}
\end{figure}

In this case the genus two Heegaard surface is unique as in Section \ref{sec2}. We regard $S^3$ as $S^3 \# S^3$ and 
we consider a standard Heegaard diagram in Figure \ref{fig:thirtytwo}. We denote by $V \cup _{\Sigma} W$ the corresponding genus two Heegaard splitting. 
We give $\partial D$ and $\partial E$ (resp. $\partial D'$ and $\partial E'$) letters $x$ and $y$ (resp. $x'$ and $y'$). 
We set $S$ to be $\Sigma \setminus (\partial D \cup \partial E)$. It is a sphere with four boundary components. 
We denote by $d^{+}$, $d^{-}$, $e^{+}$ and $e^{-}$ the boundary components of $S$ coming from $\partial D$ and $\partial E$. \\
\\
\ Let $K$ be a GOF-knot. If there is a subarc $c$ of $K$ representing $xx^{-1}$, 
we can see its corresponding arc $\alpha$ on $S$ is disjoint from $\partial D'$ and $\partial E'$ as follows: 
We assume endpoints of $\alpha$ are on $d^+$ and $e^{+}$. 
If $\alpha$ intersects $\partial D'$ or $\partial E'$ and if the intersection point nearest to $e^+$ is on $\partial E'$, 
this intersection point can be omitted by isotopy and if the intersection point nearest to $e^+$ is on $\partial D'$, 
$c$ represents word $x'yx'^{-1}$. Hence by Lemma \ref{lem1}, $K$ is a reduced form in $\{x',y'\}$. Thus except for $c$, no arcs of $K$ on $S$ intersect $\partial D'$. 
However $K$ must intersect $e^{+}$ since the word in $\{x,y\}$ represented by $K$ contains $x$ and arcs of $K$ 
which have one of the endpoints on $e^{+}$ must intersect $\partial D'$ in this situation. 
It is a contradiction. Hence we see $\alpha$ is disjoint from $\partial D'$ and $\partial E'$.
In this situation, we can replace $D$ and $E'$ with new $D$ and new $E'$ and get another standard Heegaard diagram as in Figure \ref{fig:thirtythree}. 
Note that the intersection number of $\partial ({\rm new}D) \cup \partial E$ and $K$ is less than that of $\partial D \cup \partial E$ and $K$. 
By the symmetry, if there is a subarc $c$ of $K$ representing $yy^{-1}$, 
we get another standard Heegaard diagram in which the intersection number of $\partial D \cup \partial E$ and $K$ decreases too. 
Hence we assume $K$ is a reduced form in $\{x,y\}$ on a standard Heegaard diagram.\\
\\\\
\ Changing the role of $\{\partial D,\partial E\}$ and $\{\partial D',\partial E'\}$, we assume $K$ is a reduced form in $\{x',y'\}$ in a standard Heegaard diagram.  
If there is  a subarc $c$ of $K$ representing $xx^{-1}$, we can see the corresponding arc $\alpha$ on $S$ is disjoint from $\partial D'$ and $\partial E'$ as before. 
By Lemma \ref{rem}, there must be a subarc $\bar{c}$ of $K$ representing $x^{-1}x$ and the corresponding arc on $S$ is disjoint from $\partial D'$ and $\partial E'$. 
Since $K$ is a reduced form in $\{x',y'\}$, there are no subarcs of $K$ on $S$ intersecting $\partial E'$ except for $c$ and $\bar{c}$. 
Isotoping $K$, we make $c$ and $\bar{c}$ intersect $\partial D'$ as in Figure \ref{fig:thirtyfour}. 
Since $K$ is reduced in $\{x',y'\}$, all subarcs of $K$ on $S$ except for $c$ and $\bar{c}$ are disjoint from $\partial D'$ and $\partial E'$. 
Then the other subarcs of $K$ on $S$ are like in Figure \ref{fig:thirtyfour}. 
In Figure \ref{fig:thirtyfour}, we set $n$ to be the number of arcs connecting $e^{+}$ and $d^{+}$. 
In this setting, since the number of $y$ and $y^{-1}$ in the word in $\{x,y\}$ represented by $K$ is $n$, $n$ must be even. 
Then $K$ represents the word $xx^{-1}(yx^{-1})^{\frac{n}{2}}(xy^{-1})^{\frac{n}{2}}$. It cannot be a commutator after reduction. 
Therefore $K$ does not have subarcs of type $xx^{-1}$ (and $yy^{-1}$ by symmetry). 
We assume $K$ is a reduced form in $\{x,y\}$ and $\{x',y'\}$ simultaneously.\\
\\
\ In this situation, $K$ is like in Figure \ref{fig:thirtyfive}. In this figure, the Haken sphere decomposes $K$ and its fiber into two fibered annuli in $S^3$. \\
\\
\ In $S^3$, every fibered annulus is $+1$ or $-1$ Hopf annulus. By the plumbing two fibered annuli, we obtain the trefoil or the figure eight knot. 
Hence we see every GOF-knot in $S^3$ is the trefoil or the figure eight knot. This agrees with the classical result.
\\\vspace{0.5in}

\subsection{$L(p,q)\ (|p| \geq 2)$}
\begin{figure}[htbp]
 \begin{center}
  \includegraphics[width=80mm]{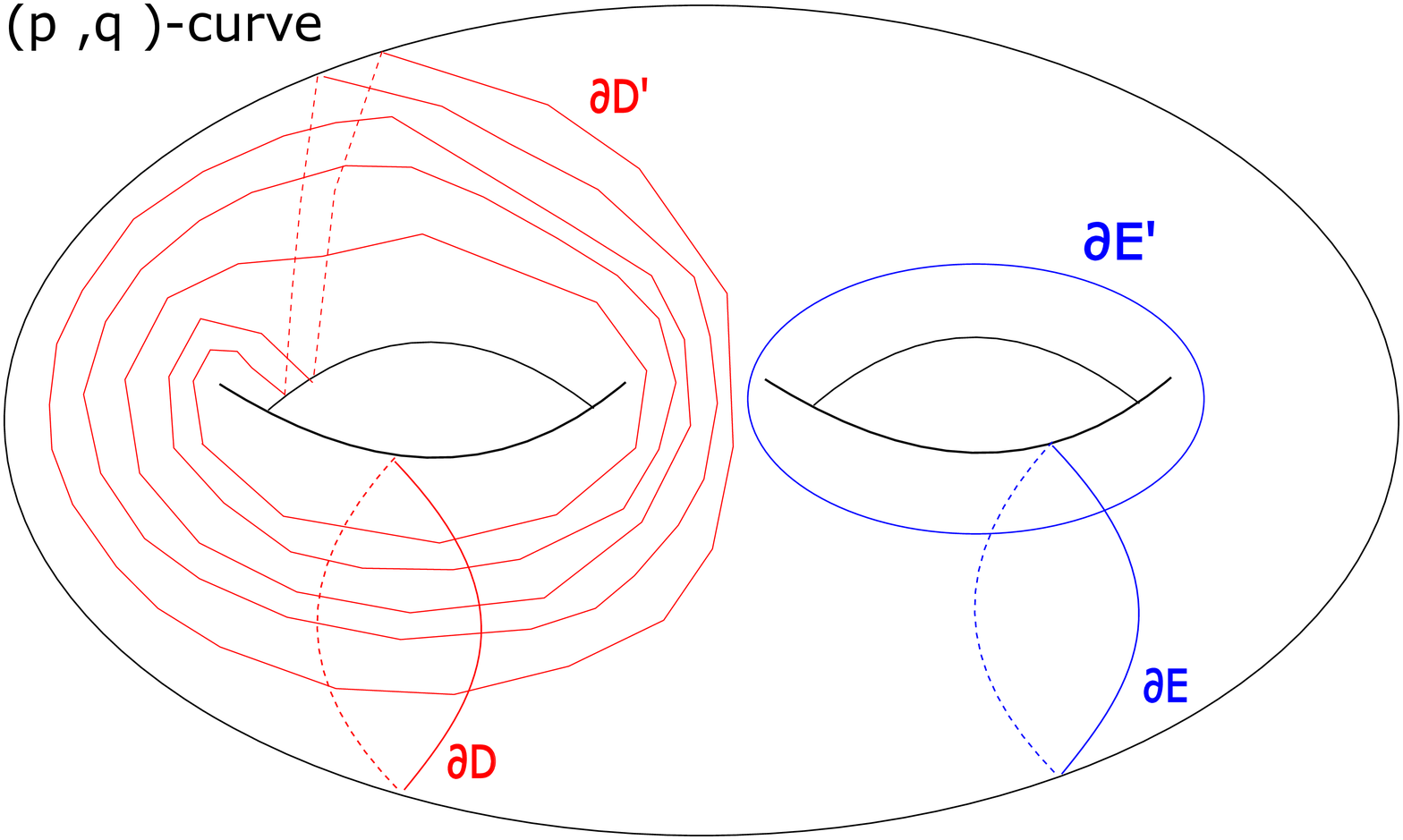}
 \end{center}
 \caption{a standard Heegaard diagram}
 \label{fig:thirtysix}
\end{figure}

\begin{figure}[htbp]
 \begin{center}
  \includegraphics[width=45mm]{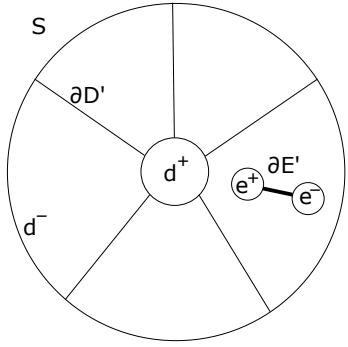}
 \end{center}
 \caption{cut $\Sigma$ along $\partial D$ and $\partial E$}
 \label{fig:thirtyseven}
\end{figure}

\begin{figure}[htbp]
 \begin{minipage}{0.5\hsize}
  \begin{center}
   \includegraphics[width=50mm]{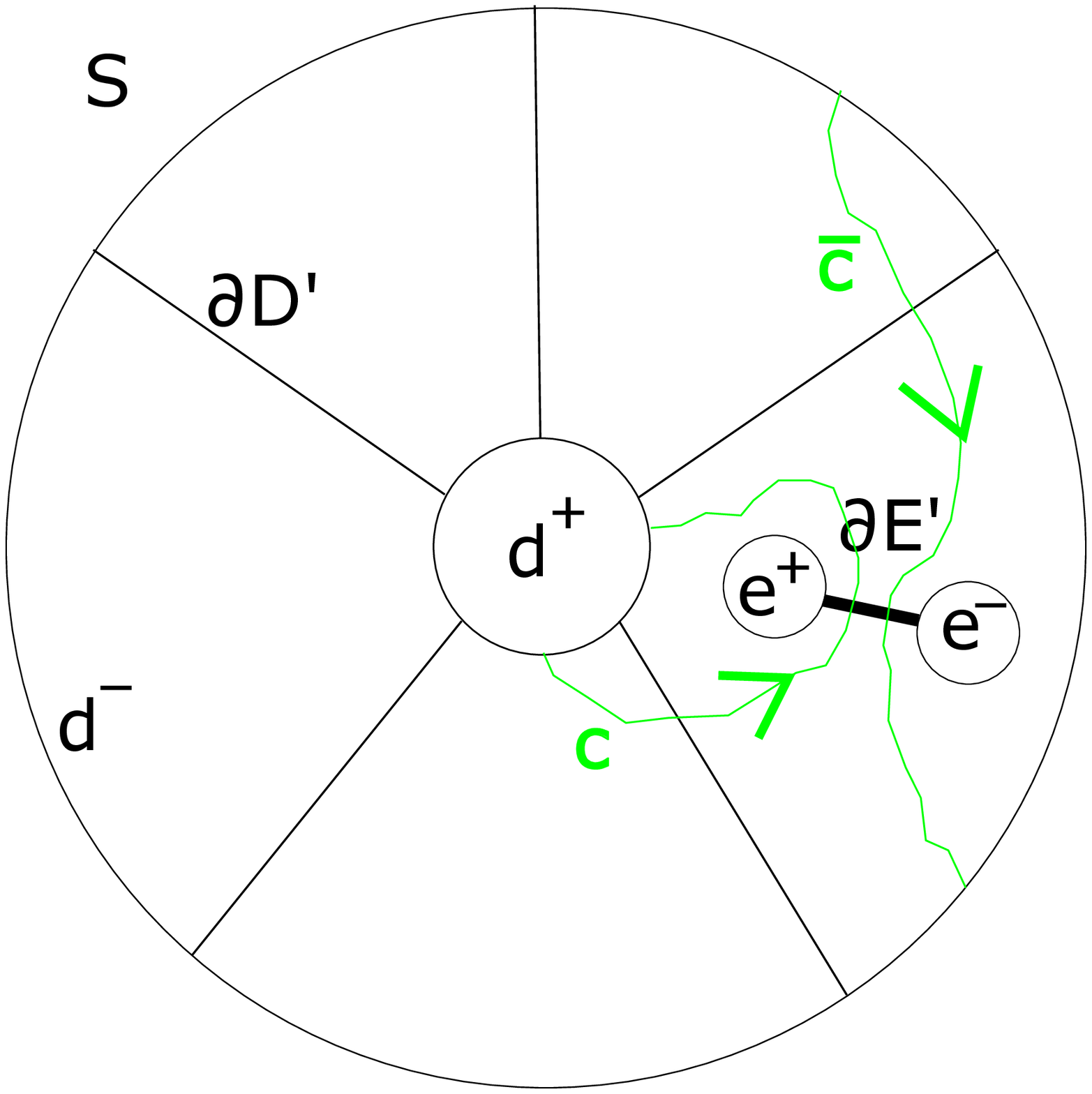}
  \end{center}
  \caption{$c$ and $\bar{c}$}
  \label{fig:thirtyeight}
 \end{minipage}
 \begin{minipage}{0.5\hsize}
  \begin{center}
   \includegraphics[width=50mm]{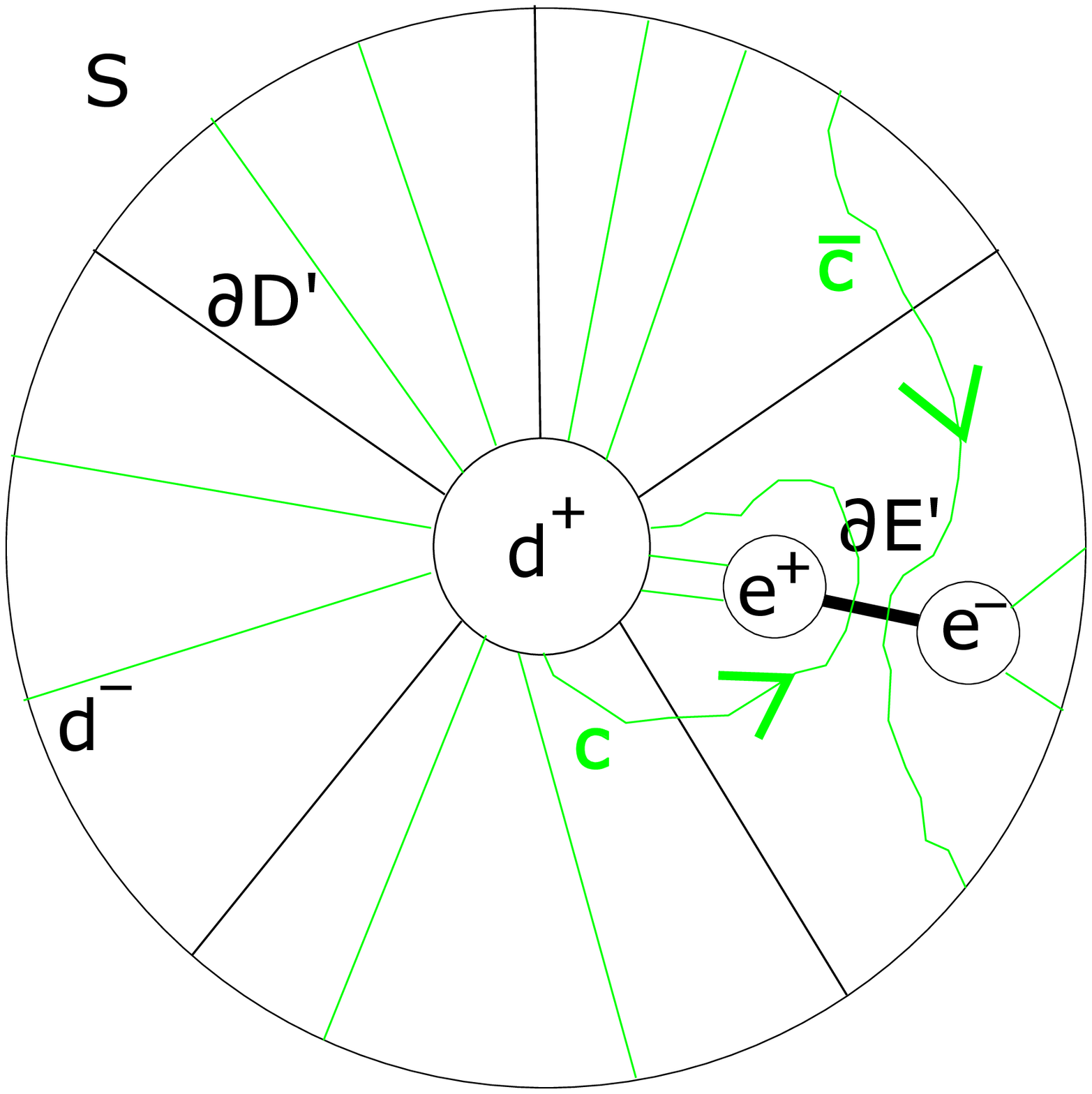}
  \end{center}
  \caption{$K$ on $S$}
  \label{fig:thirtynine}
 \end{minipage}
\end{figure}

\begin{figure}[htbp]
 \begin{center}
  \includegraphics[width=140mm]{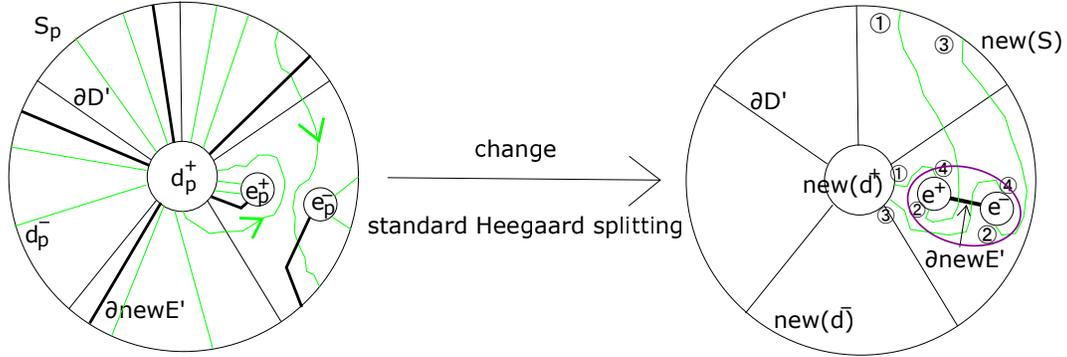}
 \end{center}
 \caption{change to another standard Heegaard splitting}
 \label{fig:forty}
\end{figure}

\begin{figure}[htbp]
 \begin{center}
  \includegraphics[width=140mm]{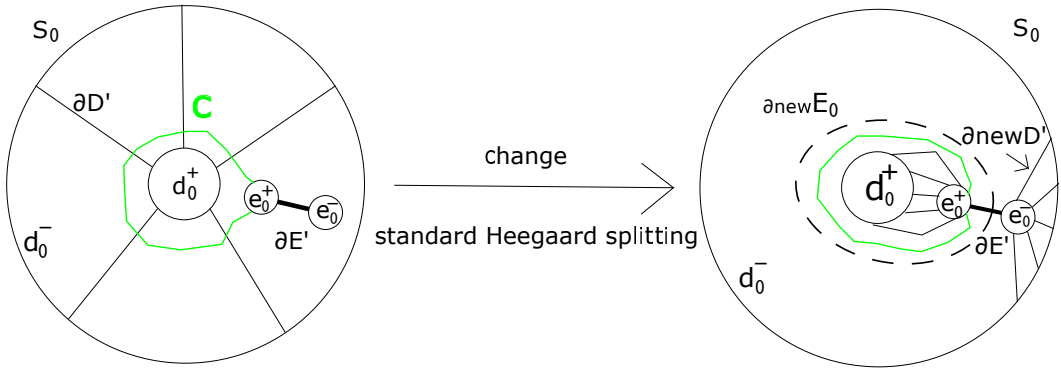}
 \end{center}
 \caption{$y_0{y_0}^{-1}$ on two standard Heegaard splittings}
 \label{fig:fortyone}
\end{figure}

\begin{figure}[htbp]
 \begin{center}
  \includegraphics[width=140mm]{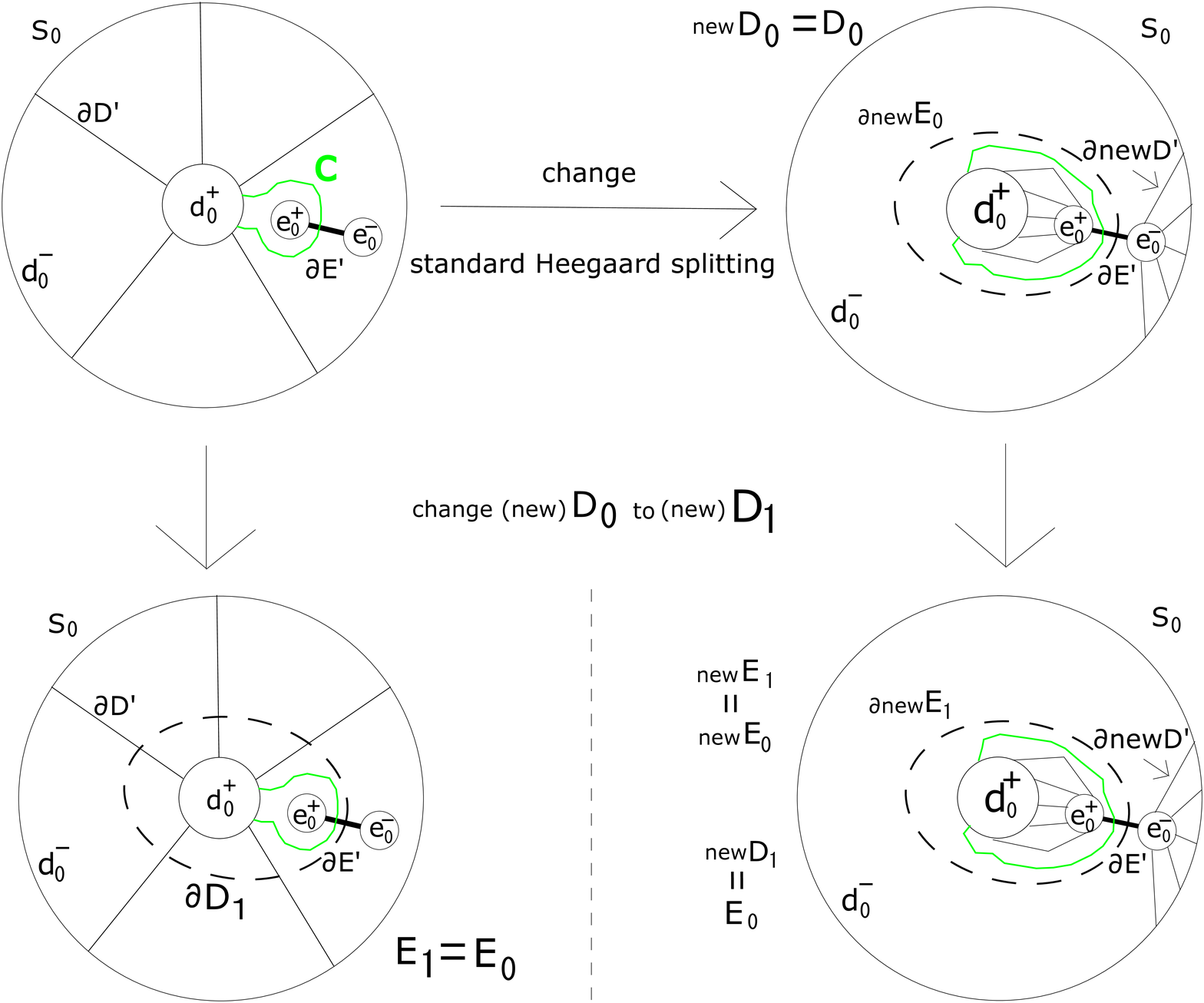}
 \end{center}
 \caption{$x_0{x_0}^{-1}$ on two standard Heegaard splittings}
 \label{fig:fortytwo}
\end{figure}

\begin{figure}[htbp]
 \begin{center}
  \includegraphics[width=120mm]{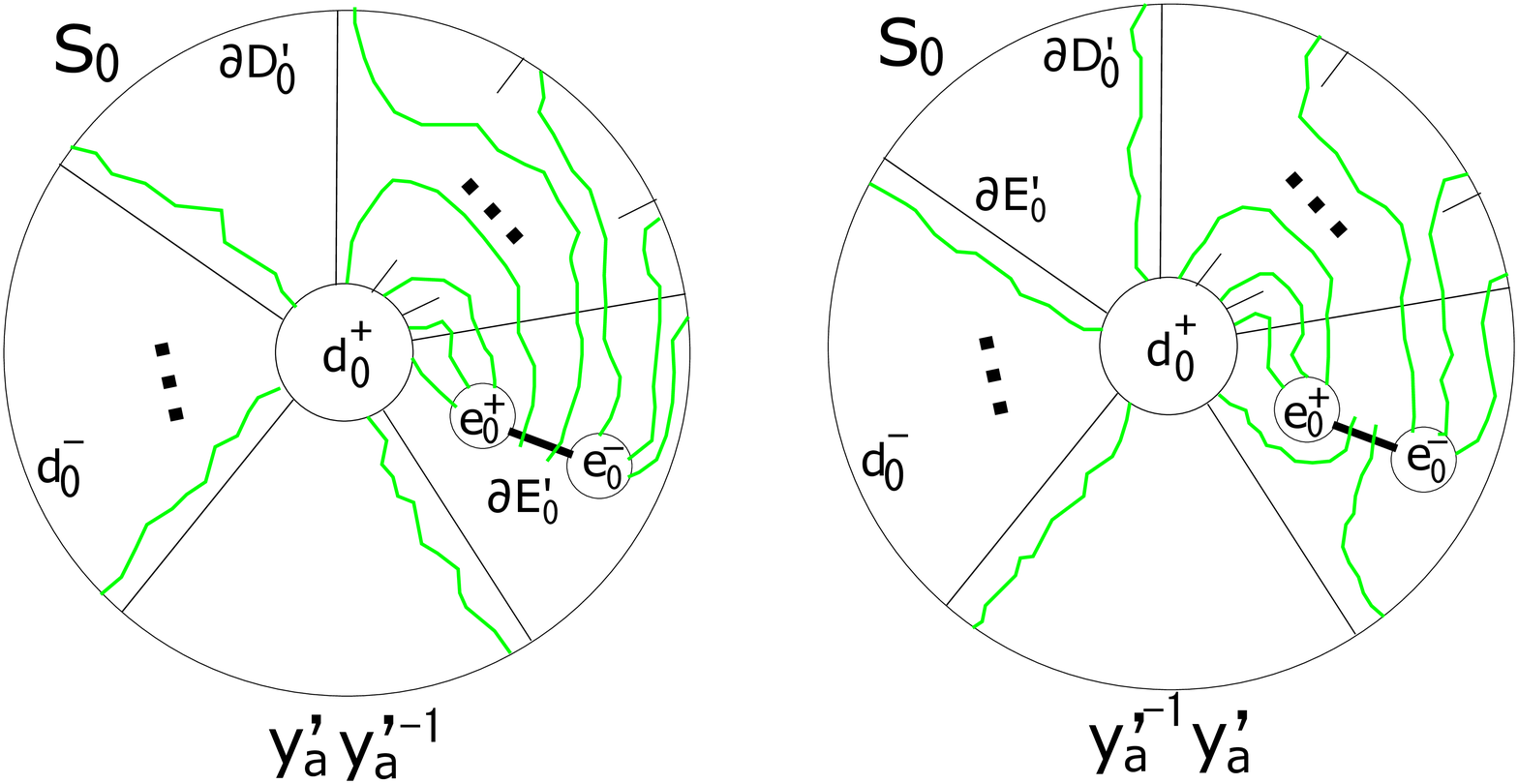}
 \end{center}
 \caption{${y'}_a{{y'}_a}^{-1}$ and ${{y'}_a}^{-1}{y'}_a$ on $S_0$}
 \label{fig:fortythree}
\end{figure}

\begin{figure}[htbp]
 \begin{center}
  \includegraphics[width=60mm]{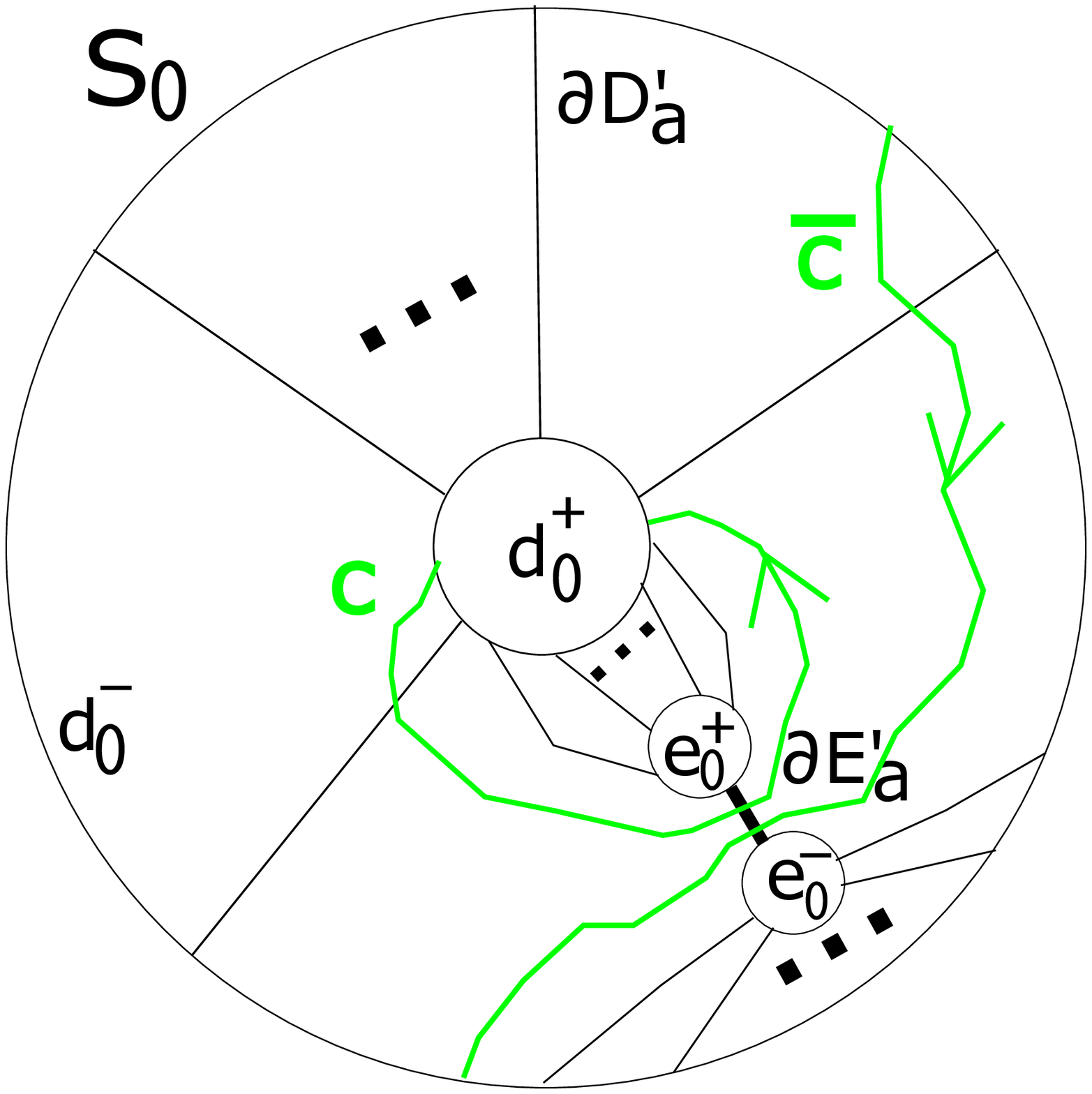}
 \end{center}
 \caption{$S_0$ with $\partial {D'}_a$ and $\partial {E'}_a$}
 \label{fig:fortyfour}
\end{figure}

\begin{figure}[htbp]
 \begin{center}
  \includegraphics[width=60mm]{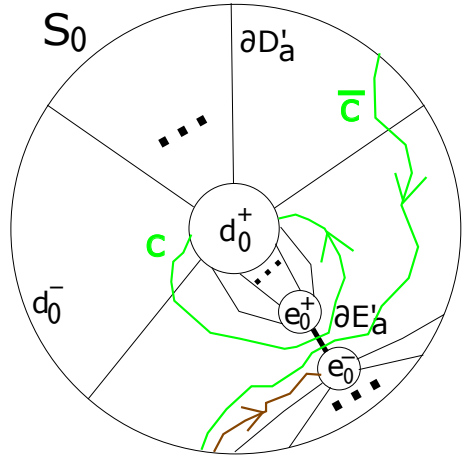}
 \end{center}
 \caption{subarc of $K$ representing ${x_0}^{-1}y_0$}
 \label{fig:fortyfive}
\end{figure}

\begin{figure}[htbp]
 \begin{center}
  \includegraphics[width=100mm]{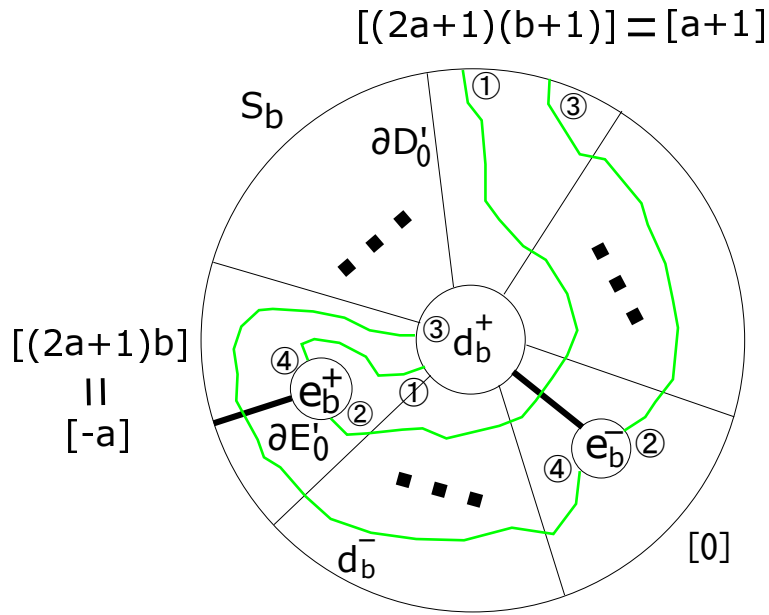}
 \end{center}
 \caption{$K$ on $S_b$ (in $L(2ab+a+b,2a+1)$)}
 \label{fig:fortysix}
\end{figure}

\begin{figure}[htbp]
 \begin{center}
  \includegraphics[width=100mm]{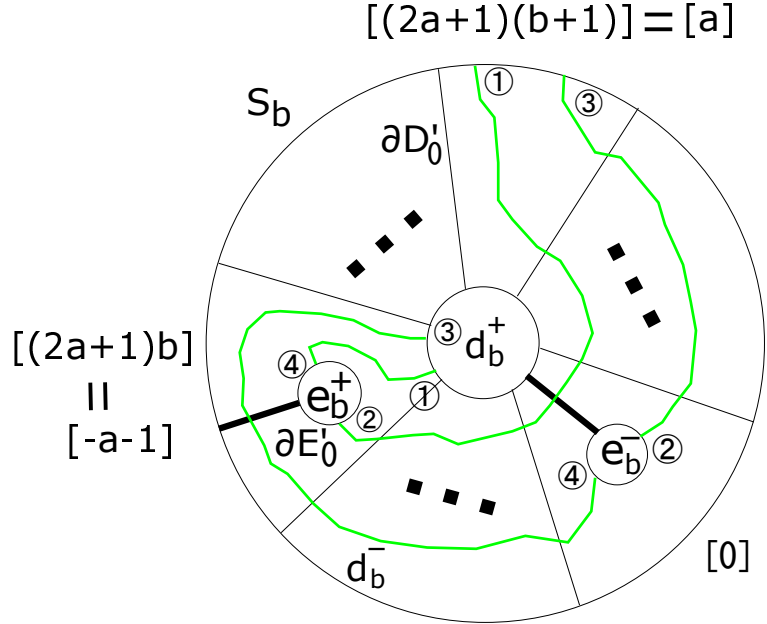}
 \end{center}
 \caption{$K$ on $S_b$ (in $L(2ab+a+b+1,2a+1)$)}
 \label{fig:fortyseven}
\end{figure}

At first, we determine the condition for $L(p,q)$ to have GOF-knots. Next, we find the positions of GOF-knots if there are. 
In this case, the Heegaard surface is unique as in Section \ref{sec2}. We consider a standard Heegaard diagram in Figure \ref{fig:thirtysix}. 
We denote by $V \cup _{\Sigma} W$ the corresponding genus two Heegaard splitting. 
We give $\partial D$ and $\partial E$ (resp. $\partial D'$ and $\partial E'$) letters $x$ and $y$ (resp. $x'$ and $y'$). 
We set $S$ to be $\Sigma \setminus (\partial D \cup \partial E)$. It is a sphere with four boundary components. 
We denote by $d^{+}$, $d^{-}$, $e^{+}$ and $e^{-}$ the boundary components of $S$ coming from $\partial D$ and $\partial E$. See Figure \ref{fig:thirtyseven}. \\
\\
\\
\ We assume that there is a GOF-knot $K$ on $\Sigma$. 
We set $S_0$ to be $S$, $d_{0}^{\pm}$ and $e_{0}^{\pm}$ to be $d^{\pm}$ and $e^{\pm}$, and $\{x_0,y_0\}$ to be $\{x,y\}$. We define $S_k$ inductively as follows. 
Suppose we have constructed up to $S_{k}$. If there are no subarc of $K$ representing the word $x_{k}x_{k}^{-1}$, we stop at getting $S_{k}$. 
If there is a subarc $c_k$ of $K$ on $S_k$ representing the word $x_kx_{k}^{-1}$ ($0 \leq k \leq p-1$), 
we can assume its corresponding arc $\alpha _k$ on $S_k$ is disjoint from $\partial D'$ and $\partial E'$. 
For, we assume $\alpha _k$ starts at $d_k^{+}$ and ends on $e_{k}^{+}$. 
If $\alpha_k$ intersects them and if the intersection point nearest to $e_{k}^{+}$ is on $\partial E'$, 
this intersection point can be omitted by isotopy and if the intersection point nearest to $e^{+}$ is on $\partial D'$, 
the subarc of $c_k$ (so of $K$) represents a word of type $x'y'x'^{-1}$ and then by Lemma \ref{lem1}, $K$ is a reduced form in $\{x',y'\}$ 
and in such case by changing the roles of $\{D,E\}$ and $\{D',E'\}$, we can assume we stop at getting $S_0$. 
Thus we assume $\alpha_k$ is disjoint from $\partial D'$ and $\partial E'$. 
We set $D_{k+1}$ to be the non-separating disk in $V$ whose boundary is the boundary of a regular neighborhood of $d_{k}^{+} \cup \alpha_k \cup e_{k}^{+}$ in $S_k$, 
and set $E_{k+1}$ to be $E_{k}$. 
We give $\partial D_{k+1}$ and $\partial E_{k+1}$ the letters $x_{k+1}$ and $y_{k+1}$. 
We set $S_{k+1}$ to be $\Sigma \setminus (\partial D_{k+1} \cup \partial E_{k+1})$. 
We denote by $d_{k+1}^{+}$, $d_{k+1}^{-}$, $e_{k+1}^{+}$ and $e_{k+1}^{-}$ the boundary components of $S_{k+1}$ coming from $\partial D_{k+1}$ and $\partial E_{k+1}$. 
As in the case of $L(p_1,q_1) \# L(p_2,q_2)$, by changing a standard Heegaard splitting, we can assume that we stop at getting $S_b$ ($0 \leq b \leq p-2$). 
(If we get $S_p$ or stop at getting $S_{p-1}$, we can assume that we stop at getting $S_0$ or $S_1$ by taking another Heegaard splitting.)
\\\\
(1) The case where we stop at getting $S_0$ and there is no subarc of $K$ representing $y_0y_0^{-1}$\\
\ \ In other words, this is the case where $K$ is a reduced form in $\{x,y\}$ or $\{x',y'\}$ on a standard Heegaard splitting. 
Changing the roles of $D$, $E$ and $D'$, $E'$ if necessary, we assume $K$ is a reduced form in $\{x',y'\}$. 
In this situation there are no subarcs of $K$ representing $yy^{-1}$. 
If there are no subarcs of $K$ representing $xx^{-1}$, $K$ is a reduced form in $\{x,y\}$ and $\{x',y'\}$ simultaneously on a standard Heegaard splitting, 
and then $L(p,q)$ must be homeomorphic to $L(p,1)$ (as in 4.2 ) and the fiber of $K$ is obtained by the plumbing of fibered annuli of $L(p,q)$ and $S^3$. 
If there is a subarc $c$ of $K$ representing $xx^{-1}$, there must be a subarc $\bar{c}$ of $K$ representing $x^{-1}x$ by Lemma \ref{rem}. 
By the reducibility of $K$ in $\{x',y'\}$, there are no subarcs of $K$ in $S$ intersecting with $\partial E'$ except for $c$ and $\bar{c}$. 
By isotoping $c$ and $\bar{c}$ so that they intersect $\partial D'$, the other subarcs of $K$ on $S$ are disjoint from $\partial D'$ and $\partial E'$. 
See Figure \ref{fig:thirtyeight}. In this figure, there cannot be subarcs of $K$ on $S$ representing $yx^{n}x^{-n}y^{-1}$ or $yx^{n}x^{-(n+m)}x^{m}y^{-1}$. 
Hence $K$ represents the word $x^{n+1}x^{-n}yx^{-n-1}x^{n}y^{-1}$ or $x^{n}x^{-n-1}yx^{-n}x^{n+1}y^{-1}$ in $\{x,y\}$. ($n$ is a natural number.) 
In particular $K$ intersects $\partial E$ twice. Drawing a picture, we see $n=p$ and $K$ represents a word $x^{-p}yx^{-l}x^{p}y^{-1}x^{l}$. 
($l$ is the minimal natural number such that $lq \equiv 1$ mod $p$.) Then $l=p\pm1$ and so $q \equiv \pm 1$ mod $p$. 
Changing the orientation if necessary, we assume $L(p,q)=L(p,1)$ and in this case $K$ is like in Figure \ref{fig:thirtynine}. 
In this figure, the operation of cancelling $x^{n+1}x^{-n}$, which changes $E'$ to a new non-separating disk, makes a new standard Heegaard splitting 
where $K$ is a reduced form in $\{x,y\}$ and $\{x',y'\}$ simultaneously. See Figure \ref{fig:forty}. 
Hence we conclude that if $K$ is a reduced form in $\{x,y\}$ or $\{x',y'\}$ in a standard Heegaard splitting, 
then $L(p,q)$ is homeomorphic to $L(p,\pm 1)$ and the fiber is constructed by the plumbing of fibered annuli in $L(p,q)$ and $S^3$. 
By changing the orientation if necessary, we assume $L(p,q)$ is $L(p,1)$ and the fibered annulus in it is the $p$-Hopf band.\\\\
\\
\\\\
(2) The case where (1) does not occur\\
\ \ In this case we stop at getting $S_b$ ($1\leq b\leq p-2$) or we stop at $S_0$ and there is a subarc of $K$ representing $y_0y_{0}^{-1}$. 
We set $D'_0$ and $E'_0$ to be $D'$ and $E'$ and give them letters ${x_0}'$ and ${y_0}'$. We set $S'_0$ to be $\Sigma \setminus (\partial D'_0 \cup \partial E'_0)$. 
We denote by ${d'_{0}}^{+}$, ${d'_{0}}^{-}$, ${e'_{0}}^{+}$ and ${e'_{0}}^{-}$ the boundary components of $S'_{0}$ coming from $\partial D'_{0}$ and $\partial E'_{0}$. 
We define ${S'}_l$ ($0\leq l\leq p$) inductively as follows. 
Suppose we have constructed up to ${S'}_{l}$. If there are no subarcs of $K$ on ${S'}_l$ representing the word ${x'_l}{x'_l}^{-1}$, we stop at getting ${S'}_l$. 
If there is a subarc $c_l$ of $K$ on $S'_l$ representing the word $x'_l{x'_l}^{-1}$ ($0 \leq l \leq p-1$), 
its corresponding arc $\alpha _l$ on $S'_l$ is disjoint from $\partial D$ and $\partial E$. 
For, we assume $\alpha _l$ starts at ${{d'}_l}^{+}$ and ends on ${{e'}_{l}}^{+}$. 
If $\alpha_l$ intersects them and if the intersection point nearest to ${e'_{l}}^{+}$ is on $\partial E$, 
this intersection point can be omitted by isotopy and if the intersection point nearest to ${e'}^{+}$ is on $\partial D$, 
the subarc of $c_l$ (so of $K$) represents a word of type $xyx^{-1}$ and then by Lemma \ref{lem1}, $K$ is a reduced form in $\{x,y\}$ 
and this contradicts our assumption. (We are in the case (2).) 
Thus we see $\alpha_l$ is disjoint from $\partial D$ and $\partial E$. 
We set $D_{l+1}$ to be the non-separating disk in $V$ whose boundary is the boundary of a regular neighborhood of ${d'_{l}}^{+} \cup \alpha_l \cup {e'_{l}}^{+}$ in $S'_l$, 
and set $E'_{l+1}$ to be $E'_{l}$. 
We give $\partial D'_{l+1}$ and $\partial E'_{l+1}$ the letters $x'_{l+1}$ and $y'_{l+1}$. 
We set $S'_{l+1}$ to be $\Sigma \setminus (\partial D'_{l+1} \cup \partial E'_{l+1})$. 
We denote by ${d'_{l+1}}^{+}$, ${d'_{l+1}}^{-}$, ${e'_{l+1}}^{+}$ and ${e'_{l+1}}^{-}$ the boundary components of $S'_{l+1}$ 
coming from $\partial D'_{l+1}$ and $\partial E'_{l+1}$. 
If this operation continues until we get ${S'}_{p}$, we can get another standard Heegaard splitting. By changing the Heegaard splitting, $E_0$ changes to new $E_0$. 
Now, in the case (2), there must be a subarc $c$ of $K$ representing $x_0{x_0}^{-1}$ or $y_0{y_0}^{-1}$. 
If $c$ represents $y_0{y_0}^{-1}$, the intersection number of $K$ with $D_0 \cup$ (new$E_0$) is less than that of $K$ with $D_0 \cup$ (old$E_0$). 
See Figure \ref{fig:fortyone}. 
If $c$ represents $x_0{x_0}^{-1}$, though the intersection number of $K$ with $D_0 \cup$ (new$E_0$) is not necessarily less than
 that of $K$ with $D_0 \cup$ (old$E_0$), $c$ represents new $x_{0}{x_0}^{-1}$ in new $S_0$ and we get new $S_1$. 
Then the intersection number of $K$ with (new$D_1) \cup$ (new$E_1$) is less than that of $K$ with (old$D_1) \cup$ (old$E_1$). See Figure \ref{fig:fortytwo}. 
Hence we see that the operation to get $S_k$ stops at getting $S_b$ ($0\leq b\leq p-2$) and 
simultaneously the operation to get ${S'}_{l}$ stops at getting ${S'}_a$ ($0\leq a\leq p-1$) (or the case (1) occurs). 
Moreover, as in the case of $L(p_1,q_1) \# L(p_2,q_2)$ we can assume the operation to get ${S'}_l$ stops at getting ${S'}_a$ ($0\leq a\leq p-2$). 
(If we stop at getting ${S'}_{p-1}$, we can assume we stop at getting ${S'}_1$ or the case (1) occurs by taking another standard Heegaard splitting.)\\
\\
\\
\ \ In this situation, if there is a subarc $c$ of $K$ representing $y_b{y_b}^{-1}$ on $S_b$, there cannot be subarcs of $K$ representing ${y'}_a{{y'}_a}^{-1}$. 
For if there is a subarc of $K$ representing ${y'}_a{{y'}_a}^{-1}$, there must be a subarc of $K$ representing ${{y'}_a}^{-1}{y'}_a$  by Lemma \ref{rem} 
and in this situation an arc representing $x_b{x_b}^{-1}$ must intersect the subarcs of $K$. See Figure \ref{fig:fortythree}.\\
\\
\ Hence we can assume $K$ is a reduced form in $\{x_b,y_b\}$ on $S_b$ or in $\{{x'}_a,{y'}_a\}$ on $S'_a$. 
We assume $K$ is a reduced form in $\{{x'}_a,{y'}_a\}$ on $S'_a$.\\
\\
\ \ Since we are in the case (2), there is a subarc $c$ of $K$ representing $x_0{x_0}^{-1}$ or $y_0{y_0}^{-1}$ on $S_0$. 
If $c$ represents $y_0{y_0}^{-1}$ in $\{x_0,y_0\}$, the letter ${x'}_a$ appears ($p-a$) times in its representing word in $\{{x'}_a,{y'}_a\}$. 
Thus $K$ cannot be a GOF-knot. 
Hence $c$ represents $x_0{x_0}^{-1}$ (and this implies $b\geq 1$). 
By Lemma \ref{rem}, there is a subarc $\bar{c}$ of $K$ representing ${x_0}^{-1}x_0$ on $S_0$. 
Each of $c$ and $\bar{c}$ intersects with $\partial {E'}_0$($=\partial {E'}_{a}$) once. 
Isotope $c$ and $\bar{c}$ so that each of them intersects with $\partial {D'}_a$ once. See Figure \ref{fig:fortyfour}. 
Since $K$ is a reduced form in $\{{x'}_a,{y'}_a\}$, the other subarcs of $K$ on $S_0$ do not intersect with $\partial {D'}_a \cup \partial {E'}_a$. \\
\\
\ \ Moreover, if $b=1$ and there is a subarc of $K$ representing $y_1{y_1}^{-1}$ on $S_1$, 
this subarc intersects with $\partial {D'}_a$ in ($p-a$) times in the same orientation, 
and if $2\leq b$ and there is a subarc of $K$ representing $y_b{y_b}^{-1}$ on $S_b$, 
this subarc intersects with $\partial {E'}_a$ at least twice in the same orientation. 
Hence we see that $K$ is a reduced form not only in $\{{x'}_a,{y'}_a\}$ but also in $\{x_b,y_b\}$. This implies $a\geq 1$. 
In particular, $K$ must intersect with $\partial E_0$($=\partial E_b$) just twice in the opposite orientation.\\\\
\ \ In Figure \ref{fig:fortyfour}, along $K$ we take some arcs on $S_0$, which are disjoint from $\partial {D'}_a \cup \partial {E'}_a$ after $c$, 
and then we take $\bar{c}$. 
Then, the terminal point of $\bar{c}$ is connected to some arc on $S_0$, 
which is disjoint from $\partial {D'}_a \cup \partial {E'}_a$ and the terminal point of this arc is connected to $\cdots$ , 
and the terminal point of this arc is connected to the starting point of $c$. 
By the symmetry of $c$ and $\bar{c}$ in $S_0$ in Figure \ref{fig:fortyfour}, 
the word in $\{x_0,y_0\}$ of the subarc of $K$, which starts at the middle point of $c$ and ends at the middle point of $\bar{c}$ is obtained by 
changing $x_0$ to ${x_0}^{-1}$ and $y_0$ to ${y_0}^{-1}$ in that of the other subarc of $K$, 
which starts at the middle point of $\bar{c}$ and ends at the middle point of $c$. 
In particular, each of these two arcs intersects $\partial E_0$ once. Since the sign of $x_0$ in each of these words (represented by these arcs) is the same, 
these words are ${x_0}^{-n}y_0{x_0}^{-m}$ and ${x_0}^{n}{y_0}^{-1}{x_0}^{m}$ ($n$ and $m$ are integers.). 
Hence the word in $\{x_0,y_0\}$ represented by $K$ is ${x_0}^{m}{x_0}^{-n}y_{0}{x_0}^{-m}{x_0}^{n}{y_0}^{-1}$. For this being commutator, $m-n=\pm1$ is necessary. 
By changing the orientations of $\partial D_0$ and $\partial E_0$ if necessary, we assume $m$ and $n$ are non-negative. 
By the construction of $S_b$, we see $(m,n)=(b+1,b)$ or $(b,b+1)$. \\\\
\ (A) $(m,n)=(b+1,b)$\\
\ \ \ By considering two components of $S_0 \setminus (\partial {D'}_a \cup \partial {E'}_a)$ 
which contain the terminal point of $c$ and the starting point of $\bar{c}$ respectively, 
and by noting that $K$ intersects $\partial D_0$ in $2b+1$ times from the middle point of $c$ to the middle point of $\bar{c}$, we see $q(2b+1)\equiv 1$ mod $p$. 
Moreover by noting that the subarc of $K$ from the middle point of $c$ to the middle point of $\bar{c}$ intersects $\partial E_0$ once, 
we see that $2b+1$ is the least natural number such that $q(2b+1)\equiv 1$ mod $p$. This implies $q'=2b+1$ in our definition of lens spaces. 
The subarc of the subarc of $K$ from the middle point of $c$ to the middle point of $\bar{c}$ corresponding to subword ${x_0}^{-1}y_0$ 
is like brown line in Figure \ref{fig:fortyfive}. 
Otherwise $K$ is a reduced form in ${S'}_l$ ($l < a$) or $qb\equiv 0$ mod $p$. From this we can see $qb=-a+up$ ($u$ is an integer). 
Then $(2b+1)q=q-2a+2up$, and since this is $1$ mod $p$ it is necessary that $q\equiv 2a+1$ mod $p$. By noting that $L(p,q)$ is homeomorphic to $L(p,q-p)$, 
we assume $q=2a+1$. Under this assumption, $K$ and $a$, $b$ are not changed. Note that if $q$ is positive, $u$ must be positive.\\
\ \ \ We will show that $u$ must be $1$. We set $v$ to be a natural number such that $wq=p+v$ ($1\leq v\leq q-1$). 
If $u\geq 2$, along the subarc of $K$ from the middle point of $c$ to the middle point of $\bar{c}$, we pay attention to $(u+1)$-st and $(u-1)$-st circuits. 
Since the subarc of $K$ from the middle point of $c$ to the middle point of $\bar{c}$ intersects $\partial E_0$ once, 
it is necessary that $a+1\leq v\leq q-1$ and $-a-v+q$ is not in $[-a,0]$. This condition is not satisfied if $q=2a+1$. 
Hence $u=1$ and $(p,q)=(2ab+a+b,2a+1)$. By drawing $K$ on $S_b$ after drawing $K$ on $S_0$ (it is almost unique.), $K$ is like in Figure \ref{fig:fortysix}. 
In this figure, $[n]$ is the residue class of $n$ mod $p$.\\
\\
\ (B) $(m,n)=(b,b+1)$\\
\ \ \ We use almost the same argument used in the case (A). As in (A), we can see that $q'=2b+1$, 
that subarcs of $K$ representing $x_0{x_0}^{-1}$, ${x_0}^{-1}$ and ${x_0}^{-1}y_0$ are like in Figure \ref{fig:fortyfive} 
and that $q(b+1)=-a+\bar{u}p$ ($\bar{u}$ is an integer). 
Then $2(b+1)q-q=-q-2a+2\bar{u}p$, and since this is $1$ mod $p$ it is necessary that $q\equiv -2a-1$ mod $p$. 
By noting that $L(p,q)$ is homeomorphic to $L(p,q-p)$, we assume $q=-2a-1$. Under this assumption, $K$ and $a$, $b$ are not changed. 
Note that if $q$ is negative, $\bar{u}$ must be negative. 
As (A), we can also see $\bar{u}=-1$. Hence $(p,q)=(2ab+a+b+1,-2a-1)$. By changing the orientation, we assume $(p,q)=(2ab+a+b+1,2a+1)$. 
By drawing $K$ on $S_b$ after drawing $K$ on $S_0$ (it is almost unique), $K$ is like in Figure \ref{fig:fortyseven}.\\
\\\\
\ Therefore we say about $L(p,q)$ and $K$ as follows: 
If $L(p,q)$ has a GOF-knot $K$ and it is on a standard genus two Heegaard surface in a reduced form in $\{x,y\}$ or $\{x',y'\}$, 
$L(p,q)$ is homeomorphic to $L(p,1)$ and a fiber of $K$ is obtained by the plumbing of the $p$-Hopf band in $L(p,1)$ and a fibered annulus in $S^3$. 
If $L(p,q)$ has a GOF-knot $K$ and it is on a standard genus two Heegaard surface in a reducible form in $\{x,y\}$ and $\{x',y'\}$, 
it is necessary that $L(p,q)$ is homeomorphic to $L(2ab+a+b,2a+1)$ or $L(2ab+a+b+1,2a+1)$ ($a$ and $b$ are positive integers.) 
and a fiber of $K$ is of the position discussed in the above argument. 
(Note that $K$ of the second case may be of the first case by using another standard Heegaard splitting.) \\\\
\ We conclude that a lens space which has a GOF-knot is homeomorphic to $L(p,1)$, $L(2ab+a+b,2a+1)$ or $L(2ab+a+b+1,2a+1)$. 
These lens spaces are classified under homeomorphisms into four classes, 
(i) $L(p,1)$ ($p\neq 4$), (ii) $L(2ab+a+b,2a+1)$ (($a$,$b$)$\neq$ ($1$,$1$)), (iii) $L(2ab+a+b+1,2a+1)$, (iv) $L(4,3)$ (it is homeomorphic to $L(4,1)$). 
On (i), there are just two GOF-knots with their fiber surfaces, one is obtained by the plumbing of the $p$-Hopf band in $L(p,1)$ and $+1$-Hopf annulus in $S^3$ 
and the other is obtained by the plumbing of the $p$-Hopf band in $L(p,1)$ and $-1$-Hopf annulus in $S^3$. 
These two can be distinguished by computing monodromies as follows:
Let $l$ be the core curve of the $p$-Hopf band in $L(p,1)$ and $l'$ be the core curve of a fibered annulus in $S^3$. 
By the plumbing we get a GOF-knot and its fiber $T$. 
We can put $l$ and $l'$ on $T$. Let $\tilde{l}$ and $\tilde{l'}$ be simple closed curves on $T$ which are obtained by moving $l$ and $l'$ along fibers. 
In $H_1(T)$, $\tilde{l} =l+l'$ and $\tilde{l'} =pl+(1+p)l'$ if a fibered annulus in $S^3$ is a $+1$-Hopf annulus, 
and $\tilde{l} =l-l'$ and $\tilde{l'} =pl+(1-p)l'$ if a fibered annulus in $S^3$ is a $-1$-Hopf annulus. 
Therefore the monodromies of these fibers are represented in $GL_{2}(\mathbb{Z})$ as $ \left(\begin{array}{ccc} 1 & p \\1 & 1+p \end{array} \right)$ 
and $ \left(\begin{array}{ccc} 1 & p \\-1 & 1-p \end{array} \right)$. They are not conjugate in $GL_{2}(\mathbb{Z})$. 
 
\begin{figure}[htbp]
 \begin{center}
  \includegraphics[width=100mm]{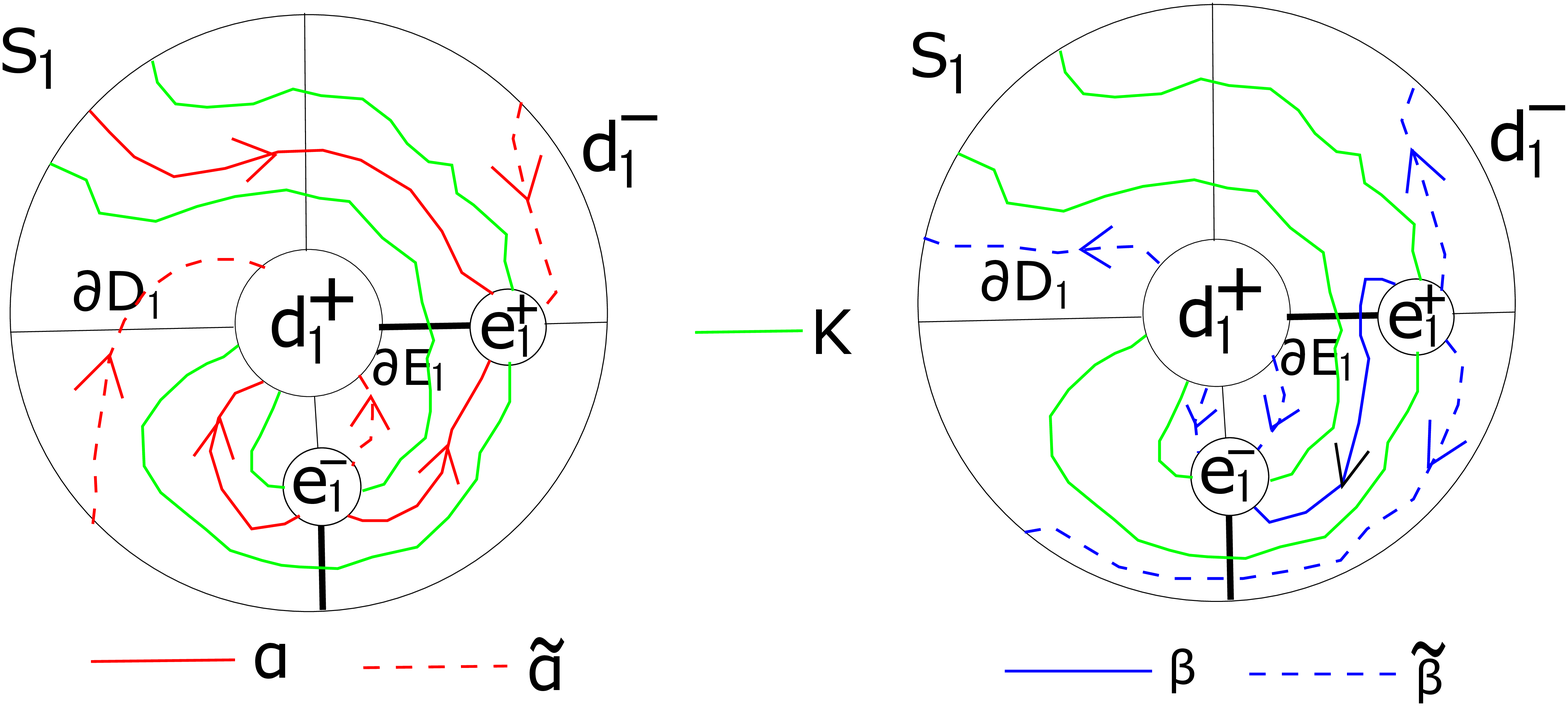}
 \end{center}
 \caption{monodromy of the third fiber in $L(4,3)$}
 \label{fig:fortyeight}
\end{figure}

On (ii) and (iii), there is just one GOF-knot, described above. 
On (iv), there are three GOF-knots, two coming from the plumbing (they can be distinguished), one coming from like (ii) ((a,b)=(1,1)). 
The third can be distinguished from the others by the monodromy. See Figure \ref{fig:fortyeight}. 
In this figure, $\alpha$ and $\beta$ are simple closed curves representing a basis of $H_1(T)$ with $T$ a fiber. 
We set $\tilde{\alpha}$ and $\tilde{\beta}$ to be simple closed curves which are obtained by moving $\alpha$ and $\beta$ along fibers. 
In $H_1(T)$, $\tilde{\alpha} = -2\alpha -3\beta$ and $\tilde{\beta} = 3\alpha +4\beta$. The monodromy of the third fiber is represented in $GL_{2}(\mathbb{Z})$ 
as $ \left(\begin{array}{ccc} -2 & 3 \\-3 & 4 \end{array} \right)$. It is not conjugate to  $ \left(\begin{array}{ccc} 1 & 4 \\1 & 1+4 \end{array} \right)$ and 
$ \left(\begin{array}{ccc} 1 & 4 \\-1 & 1-4 \end{array} \right)$ in $GL_{2}(\mathbb{Z})$. 
\\
\ This finishes a reproof of the Baker's results in \cite{3}. 
\\\\

\ GRADUATE SCHOOL OF MATHEMATICAL SCIENCES, THE UNIVERSITY OF TOKYO, 3-8-1 KOMABA, MEGURO--KU, TOKYO, 153-8914, JAPAN\\
\ \ E-mail address: \texttt{sekino@ms.u-tokyo.ac.jp}
\end{document}